\documentclass{amsart}
\usepackage[utf8]{inputenc}

\usepackage{amsmath,amsfonts, amssymb}
\usepackage{graphicx}
\usepackage{xcolor}
\usepackage{enumerate}
\usepackage{enumitem}
\usepackage{cite}
\usepackage{tabu}
\usepackage{array}
\usepackage{stmaryrd} 

\usepackage{transparent}
\usepackage{tikz}

\usepackage{fullpage}
\usepackage[active]{srcltx}

\usepackage{changepage}

\usepackage{chngcntr}

\usepackage{multibib}

\usepackage{url}
\usepackage{hyperref}
\hypersetup{
	colorlinks,
	linkcolor={red!50!black},
	citecolor={blue!50!black},
	urlcolor={blue!80!black}
}

\theoremstyle{plain}
\newtheorem{theorem}{Theorem}[section]

\newtheorem{cor}[theorem]{Corollary}

\newtheorem{lemma}[theorem]{Lemma}

\theoremstyle{definition}

\newtheorem{remark}[theorem]{Remark}

\newtheorem{construction}[theorem]{Construction}

\newcommand{\calS}{\ensuremath{{\mathcal S}}}
\newcommand{\calF}{\ensuremath{{\mathcal F}}}

\newcommand{\calT}{\ensuremath{{\mathcal T}}}

\newcommand{\comment}[1]{}

\newcommand{\bdry}{\ensuremath{\partial}}

\newcommand{\Q}{\ensuremath{\mathbb{Q}}}

\newcommand{\Z}{\ensuremath{\mathbb{Z}}}

\newcommand{\RP}{\ensuremath{\mathbb{RP}}}
\newcommand{\Sph}{\ensuremath{\mathbb{S}}}
\newcommand{\Tmax}{\ensuremath{\mathbb{T}}}

\newcommand{\mobius}{M\"{o}bius }

\definecolor{amaranth}{rgb}{0.9, 0.17, 0.31} 
\definecolor{carrotorange}{rgb}{0.93, 0.57, 0.13} 
\definecolor{citrine}{rgb}{0.89, 0.82, 0.04} 
\definecolor{dartmouthgreen}{rgb}{0.05, 0.5, 0.06} 
\definecolor{ballblue}{rgb}{0.13, 0.67, 0.8} 
\definecolor{ceruleanblue}{rgb}{0.16, 0.32, 0.75} 
\definecolor{amethyst}{rgb}{0.6, 0.4, 0.8} 
\definecolor{amber}{rgb}{1.0, 0.75, 0.0} 
\definecolor{burlywood}{rgb}{0.87, 0.72, 0.53} 

\title{Circular decompositions for hyperbolic genus one knots}

\author{Mario Eudave-Mu\~noz}

\address{Instituto de Matem\'aticas, Universidad Nacional Aut\'onoma de Mexico, Cuernavaca, Mor., MEXICO}
\email{mario@im.unam.mx}

\author{Araceli Guzm\'an-Trist\'an}
\address{Centro de Investigaci\'on en Matem\'aticas, Guanajuato, Gto., MEXICO}
\email{araceli.guzman@cimat.mx}

\author{Fabiola Manjarrez-Guti\'errez}
\address{Instituto de Matem\'aticas, Universidad Nacional Aut\'onoma de Mexico, Cuernavaca, Mor., MEXICO}
\email{fabiola.manjarrez@im.unam.mx}

\author{Enrique Ram\'{\i}rez-Losada}
\address{Centro de Investigaci\'on en Matem\'aticas, Guanajuato, Gto., MEXICO}
\email{kikis@cimat.mx}


\begin{document}
\begin{abstract}
We consider hyperbolic genus one knots having two or more non-isotopic, disjoint, genus one Seifert surfaces, with the property that the complementary regiones are handlebody pairs. We determine whether or not these knots have circular width equal to $(3)$. Moreover we give bounds on their circular widths. We also give explicit examples of such knots having three or five non-isotopic genus one Seifert surfaces and compute their circular widths. In particular this shows that many of these knots are non-almost fibered. 
\end{abstract}
\maketitle

\section{Introduction}

Let $K$ be a knot in $\Sph^3$ and $E(K)$ its exterior. A circular decomposition of $E(K)$
is defined using  collections of surfaces $\{F_i \}_1^n$ and  $\{S_i \}_1^n$, not necessarily connected, which are properly embedded in $E(K)$, such that each $F_i$ and each $S_i$ consists of a Seifert surface for $K$ and a collection (possibly empty) of closed surfaces. When the knot exterior is cut open along
the collection $\{F_i \}_1^n$ the result is a collection of disjoint submanifolds whose Heegaard surfaces (with boundary) are the $S_i$'s . A complexity $c(S_i)$ is assigned to each
$S_i$, and the circular width of $E(K)$, $cw(E(K))$, is defined. A decomposition that realizes the circular width of the knot is called circular thin. Hence when the knot exterior is in circular thin decomposition we obtain a nice sequence of Seifert surfaces which are alternately incompressible, the $F_i$'s,   and weakly incompressible, the $S_i$'s.  
The circular width of a fibered knot is defined to be zero.
An almost fibered knot is defined as a knot whose exterior possesses a
circular thin decomposition in which there is one and only one weakly incompressible Seifert
surface $S$ and one and only one incompressible Seifert surface $F$. See \cite{circulatthinpositionFMG}.

Non-fibered knots up to ten crossings are known to be almost fibered \cite{Goda}, as well as free genus one knots \cite{FreegenusoneMNR}. Examples of non-almost fibered knots can be constructed by taking the connected sum of two non-fibered knots. The existence of an infite family of genus one hyperbolic knots which are non-almost fibered was proved in \cite{MarioAraEnrique}. This family of knots contains exactly four disjoint non-parallel Seifert tori.

A genus one hyperbolic knot $K$ contains at most five non-parallel genus one Seifert surfaces, \cite{LVgenusonehypknotatmostfive}. The structure of the Kakimizu complex for these knots, $MS(K)$, has been determined \cite{LVkakimizucomplexofgenusonehypknots}. 

For a genus one hyperbolic knot $K$, let $\calT$ be the collection of isotopy classes of Seifert tori for $K$. A simplex of $MS(K)$ corresponds to a collection of mutually disjoint non-parallel Seifert tori $\Tmax= T_1 \sqcup T_2\sqcup \dots \sqcup T_N \subset E(K)$. We refer to $\Tmax$ as a \emph{simplicial collection} of $\calT$. The components are labeled consecutively as they appear around $E(K)$. The collection $\Tmax$ is maximal if its number of components $\vert \Tmax \vert=N$ is largest among all possible simplicial collections of $\calT$.  Here we show that if the result of cutting open the exterior of  $K$ along $\Tmax$ is a collection of disjoint submanifolds which are handlebody pairs, then, for many of the cases, these knots are non-almost fibered. 

When $K$ has a unique Seifert torus $T$ such that $E(K)$ cut open along $T$ is a handlebody pair, then $K$ has free genus one and by \cite{FreegenusoneMNR} it is known that $K$ is almost fibered with circular width $cw(E(K))=(3)$ or $(5)$.

The strategy to determine whether a knot has circular width $(3)$ or not, is to construct a closed manifold $\overline M$ from $E(K)$ and consider its Heegaard genus, see Construction \ref{construction}. The collection $\Tmax$ becomes a collection $\overline{\Tmax}$ of $n-1$ closed tori in $\overline M$. If we assume that $cw(E(K))=(3)$, then $\overline M$ is a genus two Heegaard manifold. Using the results of \cite{BZ} for the classification of Seifert fibered manifolds with Heegaard genus two, we are able to classify  the incompressible tori that these manifolds can contain. Then we compare $\overline{\Tmax}$ with that collection of tori to conclude that $\overline M$ is not a Seifert fibered manifold. Hence $\overline M$ has a non-trivial torus decomposition. Next we apply the classification of Heegaard genus two with non-trivial torus decomposition, \cite{KobayashiStructuresofgenustwo}, to determine whether or not $\overline M$ belongs to some of the types of such manifolds. Again we focus on the types of incompressible tori to give our conclusions. The study of such tori deeply relies in the work of \cite{LVgenusonehypknotatmostfive} and \cite{LVkakimizucomplexofgenusonehypknots}.

We can summarize our results in the following Theorems:

\begin{theorem}
Let $K\subset \Sph^3$ be a non-fibered hyperbolic genus one knot, and $\Tmax \subset \calT$ a maximal simplicial collection. Suppose that the complementary regions of $\Tmax$ are handlebody pairs.
    \begin{itemize}
        \item  If $d= dim MS(K)= 4$ or $3$, then $E(K)$ does not have circular width equal to $(3)$.
        \item  If $d= dim MS(K)= 2$, then $E(K)$ does not have circular width equal to $(3)$ unless the manifold $\overline M$ obtained by Construction \ref{construction} is of Type(iv) (see Theorem \ref{theo:kobayashi}).
        \item If $d= dim MS(K)= 1$, then $E(K)$ does not have circular width equal to $(3)$ unless the manifold $\overline M$ obtained by Construction \ref{construction} is of Type(i), Type(ii) with $M_1 \in M\ddot{o}(0,1)$, or Type (iii) with $M_1 \in D(2)$ (see Theorem \ref{theo:kobayashi}).
    \end{itemize}
\end{theorem}

\begin{theorem}\label{cor:cotasparaelanchointro}
     Let $K$ be a genus one hyperbolic knot in $\Sph^3$, $\Tmax \subset \calT$ a maximal simplicial collection, such that the complementary regions of $\Tmax$ are handlebody pairs, and $d= dim MS(K)$. 
     
     Then the circular width of $E(K)$ satisfies: 
     \begin{center}
         $(3)\leq cw(E(K))\leq (5,5,...,5)$
     \end{center}
     If all the complementary regions $\Tmax$ are simple pairs, then the circular width of $E(K)$ satisfies: 
       \begin{center}
         $(3)\leq cw(E(K))\leq (3,3,...,3)$
     \end{center}

     The right hand side of the inequalities is a $(d+1)$-tuple.
\end{theorem}

\begin{cor}
     Let $K$ be a genus one hyperbolic knot in $\Sph^3$, $\Tmax \subset \calT$ a maximal simplicial collection. If all complementary regions of $\Tmax$ are simple handlebody pairs and $(3)<cw(E(K))$, then $K$ is non-almost fibered. 
\end{cor}

We prove that genus one knots with a disjoint collection of Seifert tori such that the complementary regions are simple pairs with the power curves of every two adjacent regions intersecting in a single point are hyperbolic knots. Moreover, we describe their Kakimizu complex and their circular width is computed.

\begin{theorem} \label{theo:MaximalKakimizuSimpleCircularwidthintro}
Let $K$ be a genus one knot that contains a collection of disjoint non-isotopic Seifert tori $T_1, T_2, \dots T_n$, which decomposes $E(K)$ into simple pairs. Suppose that for every two adjacent pairs, the corresponding power curves intersect in one point. Then

\begin{enumerate}
   
    \item $K$ is a hyperbolic knot or $K$ is the trefoil and $n=1$.
    \item There is no other Seifert torus for $K$.  
    \item The Kakimizu complex of $K$, $MS(K)$, is a $(n-1)$-simplex with $n\leq 5$.
    \item If $n=1$, then $K$ is fibered or $cw(E(K))=(3)$. 
    \item If $n=2$, then $cw(E(K))=(3)$
    \item If $n=3$, $cw(E(K))=(3)$ or $(3,3)$. 
    \item If $n=4$, $cw(E(K))=(3,3)$.
    \item If $n=5$, $cw(E(K))=(3,3)$ or $(3,3,3)$.
\end{enumerate}
\end{theorem}

Lastly, we give examples of families of hyperbolic genus one knots satisfying the hypothesis of Theorem \ref{theo:MaximalKakimizuSimpleCircularwidthintro}, such that:

\begin{itemize}
\item $E(K)$ has three non-isotopic disjoint  genus one Seifert surfaces and $cw(E(K))=(3)$.
\item $E(K)$ has three  $3$ non-isotopic disjoint genus one Seifert surfaces and $cw(E(K))=(3,3)$.

\item $E(K)$ has five non-isotopic disjoint  genus one Seifert surfaces and $cw(E(K))=(3,3)$.

In this case, we conjecture that there is no knot with $cw(E(K))=(3,3,3)$ 
\end{itemize}

The examples are constructed via double branched covers of $\Sph^3$ with a $2$-component link.

   The paper is organized as follows. In Section \ref{Preliminaries} we review the circular width of knots and the known facts abour the Kakimizu complex of genus one hyperbolic knots. We also make Construction 2.5, which associates a closed manifold $\overline M$ to $E(K)$. In Sections \ref{section:SFMgenustwo} and \ref{non-trivialHeegaardgenus2} we review the incompressible tori that can exist in a Seifert fibered manifold with Heegaarde genus two, or in a manifold with a non-trivial torus decomposition and Heegaard genus two, respectively. In Section \ref{bounds} we give bounds for the circular width of hyperbolic genus one knots and prove the main results of the paper. Finally, in Section \ref{examples} we give examples of genus one knots with several Seifert tori and determine their circular width. 
   
\subsection{Acknowledgements}
EM was partially supported by grants UNAM-PAPIIT IN117423 and IN119826. MG was partially supported by grants UNAM-PAPIIT IN113323 and IN108826.  Research supported by grant SECIHTI CBF-2025-G-938.

\section{Preliminaries}\label{Preliminaries}
\subsection{Seifert tori for hyperbolic genus one knots}\label{subsection:trabajo de LuisValdez}
In this section, we state and review some of the definitions and results in \cite{LVgenusonehypknotatmostfive} and \cite{LVkakimizucomplexofgenusonehypknots}.

Let $K$ be a genus one hyperbolic knot in $\Sph^3$. Let $E(K)$ be the exterior of  $K$ in $\Sph^3$, namely the closure of the complement of an open regular neighborhood of $K$.

\begin{theorem}[Cf. \cite{LVgenusonehypknotatmostfive} Theorem 1] \label{theo:ValdezAtMost5tori}

$E(K)$ contains at most 5 mutually disjoint and non-parallel genus one Seifert surfaces.
    
\end{theorem}

In addition, in \cite{LVkakimizucomplexofgenusonehypknots} the structure of the Kakimizu complex of these knots is determined. Recall that the Kakimizu complex of a knot $K$, $MS(K)$, is the simplicial complex obtained from the isotopy classes of minimal genus Seifert surfaces for $K$, see \cite{Kaki}.

\begin{theorem}[Cf. \cite{LVkakimizucomplexofgenusonehypknots} Theorem 1]\label{theo:luisvaldez}
    If $K\subset \Sph^3$ is a genus one hyperbolic knot and  $d= dim MS(K)$ then
    \begin{itemize}
        \item $0\leq d \leq 4$,
        \item $MS(K)$ consists of at most two $d$-dimensional simplices which intersect in a common $(d-1)$-face, and is exactly one $d$-simplex if $d=0, 4$, 
        \item for each integer $1\leq d \leq 3$ there are infinitely many genus one hyperbolic knots $K\subset \Sph^3$ such that $MS(K)$ consists of two $d$-simplices
    \end{itemize}
\end{theorem}

 By Lemma 5.2 in \cite{SakumaShacktwotoriintersectatmosttwice}, the minimal intersection between two non-isotopic Seifert tori in the exterior of an atoroidal knot consists of two circles which are non-separating in each of the surfaces. 

Let $\calT$ be the collection of isotopy classes of Seifert tori for $K$, hence if two tori in $\calT$ are not disjoint, they intersect in two circles.

A simplex of $MS(K)$ corresponds to a collection of mutually disjoint non-parallel Seifert tori $\Tmax= T_1 \sqcup T_2\sqcup \dots \sqcup T_N \subset E(K)$. We refer to $\Tmax$ as a \emph{simplicial collection} of $\calT$. The components are labeled consecutively as they appear around $E(K)$. The collection $\Tmax$ is maximal if its number of components $\vert \Tmax \vert=N$ is largest among all possible simplicial collections of $\calT$. Components $T_i, T_j$ of $\Tmax$ cobound a region $R_{i,j}$ in $E(K)$ with boundary a surface of genus two that contains the longitude slope $J$ of $K$. These are called \emph{complementary regions}.

To prove his results, Valdez-Sánchez analyzes the structure of the complementary regions $R_{i,i+1}$. We will focus on the case where such complementary regions are handlebody pairs. 

A \emph{handlebody pair} consists of a genus two handlebody $H$ and a separating circle $J\subset \partial H$ which is non-trivial in $H$. 

We say that a circle in $\partial H$ is \emph{primitive} or a \emph{power} in $H$ if it represents a primitive or a power $p\geq 2$ of a non-trivial element in the fundamental group $\pi_1(H)$ of $H$, in which case the circle must be non-separating in $\partial H$.

If $(H,J)$ is a handlebody pair then the closures $T_1$ and $T_2$ of the components of $\partial H\backslash J$ are once-punctured tori with $\partial T_1=J=\partial T_2$ each of which is necessarily incompressible in $H$ and $\partial H=T_1\cup_{J}T_2$.

A handlebody pair $(H,J)$ is:

\begin{itemize}
    \item \emph{trivial} if it is homeomorphic to the pair $(T\times I, J)$ with $J$ the slope of the core of the  annulus $ \bdry T \times I$. 
\item \emph{minimal} if any once-punctured torus $T$ in $H$ with $\partial T= J$ is parallel to $T_1$
or $T_2$, so in particular any trivial pair is minimal;
\item \emph{simple} if $H = (T \times I)\cup_B V$ where $V$ is a solid torus, $B = (T \times \{0\})\cap \partial V$ is a closed annular
neighborhood of a non-separating circle $\omega\times \{0\}$ in $T \times \{0\}$, and $B$ runs $p \geq 2$ times around $V$. The
separating circle $J$ corresponds to $\partial T \times \{1/2\}$. Let $T_1 = T\times \{1\}$ and $T_2=(T\times \{ 0 \} \backslash B )\cup \bdry V \backslash B$. The core of the annulus $\bdry V \backslash B \subset T_2$ is coannular in $H$
to $\omega \times \{1\} \subset T_1$. The \emph{index} of the simple pair $(H, J)$ is defined as the integer $p\geq 2$. The solid torus $V$, or a solid torus $V'$ contained in and concentric to $V$, is called a \emph{companion torus}. Note that if $V'$ is a companion torus contained in the interior of $H$, there are disjoint annuli $A_0$ and $A_1$, contained in $H$, such that a boundary component of $A_i$ lie in $T_1$ and the other in $T_2$, for $i=0,1$. Note that the boundary component of $A_i$ lying in $T_1$ is a power in $H$, which in fact represents a power $p\geq 2$ of a non-trivial element of $\pi_1(H)$. 

An alternative description of a simple pair is as follows. Consider  a copy of $\omega \times \{0\} $ in the interior of $T\times I$, say $\omega \times \{ 1/2 \}$. Perform $p/q$-Dehn surgery on $\omega \times \{ 1/2 \}$ (taking a natural meridian and a copy of $\omega$ as a longitude). Denote the simple pair by $(H(p/q),J)$. Note that this pair has index $p$.

\item \emph{primitive} if $(H,J)$ is non-trivial and there is a non-separating annulus $A \subset H$ with each
boundary component $\partial_1 A \subset T_1$ and $\partial_2 A \subset T_2$ a primitive circle in $H$. 
\end{itemize}

$H[J]$ denotes the $2$-handle addition to $H$ along $J$. Observe that  $T_1$ and $T_2$ after the $2$-handle addition become closed tori $\hat{T}_1$ and $\hat{T}_2$, respectively. 

\begin{lemma}[Cf. \cite{LVgenusonehypknotatmostfive} Lemma $3.7(4)$]\label{lemma:trivialpairifandonlyifSeiferttoriparallel}
    Let $(H,J)$ be a pair with $\partial H=T_1\cup_J T_2$ and $T\subset H$ any once-punctured
torus with $\partial T=J$. Then the pair is trivial if and only if $H[J]=\hat{T}_1\times I$.
\end{lemma}

\begin{lemma}[Cf. \cite{LVkakimizucomplexofgenusonehypknots}Lemma $2.3.11$ ]\label{lemma:twohandleadditionforpairs}
    If $(H,J)$ is a minimal non-trivial handlebody pair then
    
    $H[J]=\begin{cases}
        \mathbb{A}(q/p) & \text{ if } (H,J) \text{ is a simple pair of index $p\geq 2$} \\
        \text{toroidal} & \text{ if } (H,J) \text{ is a primitive pair} \\
        \text{hyperbolic} & \text{ otherwise }
    \end{cases}$
\end{lemma}

where  $\mathbb{A}(q/p)$ denotes a Seifert fibered space with base space an annulus and one exceptional fiber of order $p$. Furthermore, the coefficient of the exceptional fiber is $q/p$ if $(H,J)=(H(p/q),J)$. 

\subsection{Circular width for knots in $\Sph^3$.}\label{subsection:circularwidth}

In this section, we give an overview about circular thin decomposition for knots. The reader is referred to \cite{circulatthinpositionFMG} for further exploration of the concepts. 

Let $E(K)$ be the exterior of a knot $K$ in $\Sph^3$.

Let $F$ and $F'$ be two disjoint properly embedded surfaces in $E(K)$, not necessarily connected, such that $F$ and $F'$ contain a Seifert surface for $K$. Let $R_{F,F'}$ a region in $E(K)$ bounded by $F$ and $F'$.
If there is a properly embedded surface $S$ in $E(K)$, disjoint from $F$ and  $F'$, not necessarily connected that contains a Seifert surface for $K$, and compression bodies $A$ and $B$ such that:
\begin{itemize}
   
        \item $R_{F,F'}=A\cup_SB$
        \item $\bdry_- A=F$, $\bdry_+A=S$, $\bdry_vA= \bdry F \times I$
        \item $\bdry_-B=F'$, $\bdry_+B=S$, $\bdry_v B=\bdry F'\times I$.
    \end{itemize}

We say that $S$ \emph{determines a decomposition for } $R_{F,F'}$, and that \emph{$(A,B;S)$ is a decomposition for $R_{F,F'}$}. 

The \emph{complexity of such decomposition} is $c(S)=  \sum (1-\chi(\overline S))$ where $\overline{S}\neq \Sph^2$ is the surface $S$ with its boundary capped with a disk, and the summatory runs over all the components of $\overline{S}$.

Let $\calF=\{F_1, F_2 \dots F_n \}$ and $ \calS=\{S_1, S_2, \dots S_n \}$ be collections of properly embedded disjoint surfaces in $E(K)$, not necessarily connected, such that each $F_i$ and each $S_i$ contain a Seifert surface for $K$. If  $S_i$ determines a decomposition of $R_{F_i,F_{i+1}}$ for $i=1,2,...,n-1$, then the triple $(E(K), \calF, \calS)$  is called a \emph{circular decomposition of $E(K)$}.

The surfaces in the collection $\calF$ are called \emph{thin surfaces}, while those contained in $\calS$ are called \emph{thick surfaces}.

The \emph{circular width of $(E(K), \calF, \calS)$} is an $n$-tuple given by the collection of $\{ c(S_i)\}_{i=1}^n$ ordered non-increasingly. The circular width of $(E(K), \calF, \calS)$  is denoted by $cw(E(K), \calF, \calS)$. 

Given two circular decompositions for $E(K)$ and its circular widths, we compare them using the lexicographic order. The \emph{circular width of $E(K)$} is defined as the minimal $cw(E(K), \calF, \calS)$ over all possible circular decompositions. The circular width of $E(K)$ is denoted by $cw(E(K))$. 

A circular decomposition realizing $cw(E(K))$ is a \emph{circular thin decomposition of $E(K)$}.

If a knot $K$ is fibered, we define the circular width of $K$, $cw(E(K))$, to be equal to zero. 

A knot $K$ is \emph{almost fibered} if it has a circular thin decomposition consisting of only one thin surface and only one thick surface, in other words, its circular width is a 1-tuple. 

A circular thin decomposition $(E(K), \calF, \calS)$ gives a nice sequence of Seifert
surfaces which are alternately incompressible and weakly incompressible, see \cite{circulatthinpositionFMG} Theorem 3.2.

\subsection{From $E(K)$ to a closed manifold.}
Let $K\subset \Sph^3$ be a hyperbolic genus one knot, and $\Tmax \subset \calT$ a maximal simplicial collection. Suppose that the complementary regions of $\Tmax$ are handlebody pairs.

\begin{construction}\label{construction}\ \
\\
\begin{enumerate}
    \item Perform $0$-surgery on E(K) to obtain a closed manifold $M'$ with a collection of closed incompressible tori $\Tmax'$ (see \cite{GabaiFoliations3}) which correspond to the images after the $0$-surgery of $\Tmax$.  
    \item Let $M''$ be the manifold obtained by cutting $M'$ along $T'$, where $T' \in \Tmax'$ is the image of some $T$ after $0$-surgery. The manifold $M''$ has two boundary components consisting of copies $T_1'', T_2''$ of $T'$.

\item Fill in each boundary component $T_1''$ and $T_2''$ of $M''$ with solid tori $V_1$ and $V_2$ to obtain a closed manifold $\overline M$. 

The filling can be chosen in such a way that after the filling, the collection $\Tmax'\setminus T'$ becomes a family $\overline \Tmax$ of $n-1$ incompressible closed tori for $n=2,3,4,5$ (\cite{CGLSdenhsurgeryofknots}). 
\end{enumerate}    

\end{construction}

\begin{remark}\label{remark:observacionesdespuesdelaconstruccion}
    Since the complementary regions of $\Tmax \subset E(K)$ are handlebody pairs and the tori in $\Tmax$ are disjoint and non-parallel, then:

    \begin{itemize}
        \item Lemma \ref{lemma:trivialpairifandonlyifSeiferttoriparallel} implies that the collection of closed tori $\overline \Tmax$ in $\overline M$ are non-parallel.
        \item Lemma \ref{lemma:twohandleadditionforpairs} implies that $\overline M$ is a graph manifold if each complementary region of $\Tmax \subset E(K)$ is a simple pair. 
        \item Let $\overline \calT$ denote the result of the Seifert tori $\calT$ after the construction. The closed tori in $\overline \calT$ are separating in $\overline M$, they are either disjoint or at most two of them intersect minimally in two circles which are non-separating in each of the surfaces.

    \end{itemize}
\end{remark}

\begin{remark}\label{remark:aplicaciondecontruccionagenero2}
If additionally we assume that $cw(E(K))=(3)$, then $E(K)$ has a circular thin decomposition $(E(K), T, S)$ where $T\in \Tmax$  and $S$ is a genus two Seifert surface. Both surfaces have parallel boundaries on $\bdry E(K)$, which are $0$-slopes. Moreover $E(K)\setminus(T \cup S)$ is made of two compression bodies $A$ and $B$. The closed manifold $\overline M$ obtained in Construction \ref{construction} is of Heegaard genus two. In step $(1)$ a closed genus two surface $S'$ arises from $S$. Observe that $S'$ separates $M''$ into two compression bodies $A''$ and $B''$, which arise from the compression bodies $A, B$ after capping off their vertical boundaries with a piece of the surgered solid torus. Note that $\bdry_+A''= \bdry_+ B''= S'$ and $\bdry_-A''= T_1''$ and $\bdry_- B''= T_2''$. Finally, in step $(3)$ the surface $\overline S$, the image of $S'$, is a genus two Heegaard surface for $\overline M$, since it separates $\overline M$ into two handlebodies $\overline A= A'' \cup_{T_1''} V_1$ and $\overline B= B''\cup_{T_2''} V_2$.
    
\end{remark}

\section{Incompressible tori in Seifert fibered manifolds with Heegaard genus two}\label{section:SFMgenustwo}
Boileau and Zieschang in \cite{BZ} calculate the Heegaard genus of almost all closed orientable Seifert fibered manifolds. 

We use their results to derive properties that will be helpful in proving that a manifold obtained as in Remark \ref{remark:aplicaciondecontruccionagenero2} is not a Seifert fibered manifold. To make it easier for the reader to understand the results derived from their work, we include the original statements of the main theorems. 

Following their notation, let $S(\pm g; e_0; \beta_1/\alpha_1, \dots \beta_m/\alpha_m)$ denote an orientable Seifert fibered manifold with a closed orientable basis of genus $g$, or with a closed non-orientable basis of genus $-g$. The number $e_0 \in \Q$ is the rational Euler number of the Seifert fibration and $ \beta_i/\alpha_i \in (\Q/\Z)^*$ are the exceptional fibers. The Euler number of the Seifert fibering is $e= e_0 + \sum \frac{\beta_i}{\alpha_i} \in \Z$. The rank of the fundamental group of $M$ is denoted by $r(M)$, and the Heegaard genus of $M$ by $h(M)$.

\begin{theorem}[Cf.\cite{BZ}Theorem 1.1]\label{theo:BZorientablebasis}
Let $M=S(g; e_0; \beta_1/\alpha_1, \dots \beta_m/\alpha_m)$ with $1<\alpha_1\leq \dots \leq \alpha_m$ and $gcd(\alpha_i,\beta_i)=1$, $1\leq i\leq m$.

\begin{enumerate}
    \item If $M=S(g;e_0; \frac{1}{2},\dots \frac{1}{2}, \frac{\beta_m}{2\lambda+1})$, with $\lambda>0$, $m$ even and $m\geq 4$, then $m-2 =r(M)\leq h(M)\leq m-1$.
    If in addition, $e_0=\pm \frac{1}{2(2\lambda+1)}$, then  $r(M)=h(M)=m-2$.
    \item If $g>0$ and $m=2$, or if $m\geq3$ and if $M$ does not belong to the case (1), then  $r(M)=h(M)=2g+m-1$.
    \item If $m=1$, then $r(M)=h(M)=2g$ if $e_0=\pm \frac{1}{\alpha}$; otherwise $r(M)=h(M)=2g+1$.
    If $m=0$, then $r(M)=h(M)=2g$ if $e_0=\pm 1$; otherwise $r(M)=h(M)=2g+1$.
\end{enumerate}
\end{theorem}

\begin{theorem}[Cf.\cite{BZ}Theorem 2.1]\label{theo:BZnonorientablebasis}

Let $M=S(-g; e_0; \beta_1/\alpha_1, \dots \beta_m/\alpha_m)$ be a Seifert manifold with non-orientable basis and with $1<\alpha_1\leq \dots \leq \alpha_m$ and $gcd(\alpha_i,\beta_i)=1$, $1\leq i\leq m$. 
\begin{enumerate}
    \item If $m \geq 2$, then $r(M)=h(M)=g+m-1$.
    \item If $m=1$, then $r(M)=h(M)=g$ if $e_0=\frac{\pm 1}{\alpha_1}$ and $r(M)=h(M)=g+1$ otherwise.
    \item If $m=0$, then $r(M)=h(M)=g$ if $e_0= \pm 1$ and $r(M)=h(M)=g+1$ otherwise.
\end{enumerate}
\end{theorem}

We are interested in those orientable Seifert manifolds with Heegaard genus two. Hence, we have the following corollary:

\begin{cor}\label{coro:genustwoSeifertmanifolds}
    Let $M$ be an orientable  Seifert manifold with Heegaard genus two. Then $M$ is one of the following, in each case, the basis and the exceptional fibers are described:  
    \begin{description}
        \item[$\Sph^2(4)$] A sphere with four exceptional fibers $(\frac{1}{2}, \frac{1}{2}, \frac{1}{2}, \frac{\beta}{2\lambda+1})$ and $e_0=\pm \frac{1}{2(2\lambda+1)}$
        \item [$\Sph^2(3)$] A sphere with three exceptional fibers. 
        \item [$T(1)$]A torus with one exceptional fiber and $e_0= \pm 1/\alpha_1$.
        \item[$T(0)$] A torus with zero exceptional fibers and $e_0=\pm 1$.
        \item [$\RP^2(n)$] A projective plane with $n=0,1,2$ exceptional fibers.
        \item [$K(0)$]A Klein bottle with zero exceptional fibers and $e_0=\pm 1$.
        \item [$K(1)$]A Klein bottle with one exceptional fiber and $e_0= \pm 1/\alpha_1$
    \end{description}
        
\end{cor}

Next, we investigate the incompressible tori in each of the cases of Corollary \ref{coro:genustwoSeifertmanifolds}. Moreover, we determine the number of such tori and if they are separating or not. 

\begin{figure}
        \centering
        \includegraphics[width=0.3\linewidth]{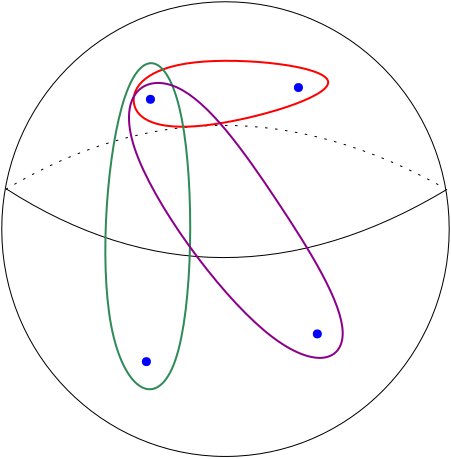}
        \caption{$\mathbb{S}^2(4)$. The three circles lift to three separating incompressible tori, every two of them intersect in two essential curves.}
        \label{fig:S(4)}
    \end{figure}

\begin{figure}
        \centering
        \includegraphics[width=0.4\linewidth]{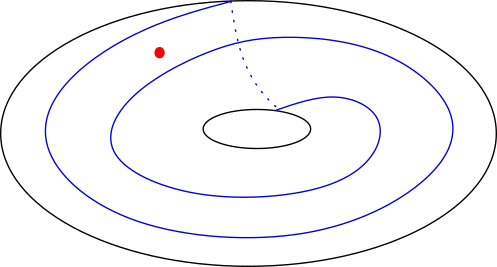}
        \caption{$T(n)$ with $n=0,1$. There are infinitely many non-separating, incompressible tori.}
        \label{fig:T(1)}
    \end{figure}

\begin{figure}
         \centering
         \includegraphics[width=0.3\linewidth]{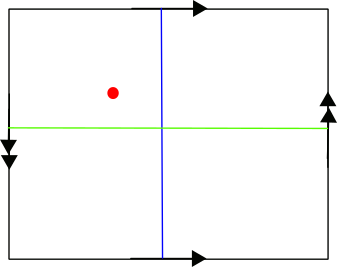}
         \caption{$K(n)$ with $n=0, 1$. The blue circle lifts to a separating incompressible torus, while the green one lifts to a Klein bottle, whose regular neighborhood has an incompressible non-separating torus as boundary. These two incompressible tori intersect in two essential curves.}
         \label{fig: Klein}
     \end{figure}

\begin{lemma}\label{lem:typeoftoriinspheretorusklein}
    The cases $\Sph^2(4), T(n), K(n)$ with $n=0,1$ do not contain horizontal surfaces. Each of them contains the following vertical tori:
    \begin{itemize}
        \item $\Sph^2(4)$ has three separating incompressible tori, every two of them intersect in two essential curves.
        \item $T(n)$ with $n=0, 1$ contains infinitely many non-separating incompressible tori.
        \item $K(n)$ with $n=0,1$ has one non-separating and one separating incompressible tori that intersect in two curves.
    \end{itemize}
    
\end{lemma}
\begin{proof}
    For all the cases in the statement, the rational Euler number is not an even integer. Hence, they do not contain horizontal surfaces. The existence and properties of the vertical tori, in each case, follow from examining Figures \ref{fig:S(4)},\ref{fig:T(1)} and \ref{fig: Klein}.
\end{proof}

\begin{lemma}\label{lem:torisphere(3)}
     If  $\Sph^2(3)$ is one of the following:
    \begin{itemize}
        \item $M_1=S(0; 0; 2/3, -1/3, -1/3)$
        \item $M_2=S(0; 0; 1/2, -1/4, -1/4)$
        \item $M_3=S(0; 0; 1/2, -1/3, -1/6)$
    \end{itemize}
    Then $\Sph^2(3)$ contains a unique horizontal surface, which is a separating torus.
    Furthermore, in any case, $\Sph^2(3)$ does not contain any incompressible vertical torus. 
\end{lemma}
\begin{proof}
    For the first part, see Section 2.2 in \cite{Hatcher}. Any incompressible vertical torus arises as the lift of an essential curve, or a curve not bounding a single point that lifts to an exceptional fiber. In $\Sph^2(3)$ such kinds of curves do not exist. Hence, $\Sph^2(3)$ does not contain any incompressible vertical torus.
\end{proof}

\begin{lemma}\label{lem:notorirp2(0,1)}
 $\RP^2(n)$ with $n=0,1$ does not contain incompressible tori.
\end{lemma}
\begin{proof}
The Seifert fibered manifold $\RP^2(n)$ with $n=0,1$ is  $M= \RP^3 \sharp \RP^3$, see \cite{Jaco}. If $M$ contains an incompressible torus, by an innermost argument and the incompressibility of the torus, it must be contained in a summand $\RP^3$, but this does not contain an incompressible torus. 
\end{proof}

\begin{figure}
         \centering
         \includegraphics[width=0.3\linewidth]{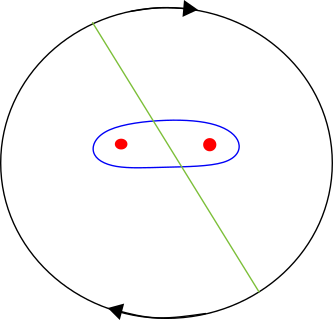}
         \caption{$\RP^2(2)$. The blue circle lifts to a separating incompressible torus. The green circle lifts to a Klein bottle whose regular neighborhood has an incompressible torus as boundary. These two incompressible tori intersect in four essential curves.}
         \label{fig:RP(2)}
     \end{figure}

 \begin{lemma}\label{lem:torirp2(2)}
      $\RP^2(2)$ contains two separating incompressible tori which intersect in four essential curves. If $\RP^2(2)$ is the manifold $M=S(-1;0; 1/2, -1/2)$, then it contains one horizontal separating torus. 
 \end{lemma}
\begin{proof}
     The vertical tori arise from the curves in the base space $\RP^2$ depicted in Figure \ref{fig:RP(2)}. For the existence of the horizontal separating torus in the case that $M=S(-1;0; 1/2, -1/2)$, see \cite{Hatcher} section 2.2.
\end{proof}

\begin{cor}\label{coro:theclosedmanifoldforextKisnotSF} 
Let $\overline M$ be as in Remark \ref{remark:aplicaciondecontruccionagenero2}, then $\overline M$ is not a Seifert fibered manifold. In other words, $\overline M$ has a non-trivial torus decomposition. 
    \end{cor}

\begin{proof}
If the manifold $\overline M$ were a Seifert fibered manifold, it would be one of the ones listed in Corollary \ref{coro:genustwoSeifertmanifolds}, since $\overline M$ has a Heegaard decomposition of genus two. Nevertheless, none of the incompressible tori in $\overline \calT$ contained in $\overline M$ are as those described in Lemmas \ref{lem:typeoftoriinspheretorusklein}-\ref{lem:torirp2(2)}. Therefore, $\overline M$ is not a Seifert fibered manifold.
\end{proof}

\section{Incompressible tori in manifolds with non-trivial torus decomposition of Heegaard genus two}\label{non-trivialHeegaardgenus2}

The Haken manifolds with non-trivial torus decomposition, that is a non-empty JSJ-decomposition,  which possesses a genus two Heegaard splitting were determined by Kobayashi in \cite{KobayashiStructuresofgenustwo}.

\begin{theorem}[Cf.\cite{KobayashiStructuresofgenustwo}]\label{theo:kobayashi}
Let $M$ be a closed, connected Haken manifold with a Heegaard splitting of genus two. If $M$ has a non-trivial torus decomposition, then either:
\begin{itemize}
    \item Type (i): $M$ is obtained by $M_1 \in D(2)$ and $M_2 \in L_K$ by identifying their boundaries, where the regular fiber of $M_1$ is identified with the meridian loop of $M_2$.
    \item Type (ii): $M$ is obtained from $M_1 \in M\ddot{o}(n)$ with $n=0,1$ or $2$ and $M_2 \in M_K$ by identifying their boundaries, where the regular fiber of $M_1$ is identified with the meridian loop of $M_2$.
    \item Type (iii): $M$ is obtained from $M_1 \in D(n)$ with $n=2$ or $3$ and $M_2 \in M_K$ by identifying their boundaries, where the regular fiber of $M_1$ is identified with the meridian loop of $M_2$.
    \item Type (iv): $M$ is obtained from $M_1, M_2 \in D(2)$  and $M_3 \in M_L$ by identifying their boundaries where the regular fiber of $M_i$, for $i=1,2$ is identified with the meridian loop of $M_3$.
    \item Type (v): $M$ is obtained from $M_1 \in A(n)$ with $n=0, 1$ or $2$ and $M_2 \in M_L$ by identifying their boundaries where the regular fiber of $M_1$ is identified with the meridian loop of $M_2$.
\end{itemize}
Conversely, if a $3$-manifold has a decomposition of Type (i)-(v), then it has a Heegaard splitting of genus two.
    
\end{theorem}

$D(n), A(n), M\ddot{o}(n)$ is the collection of Seifert fibered manifolds with base space a disk, an annulus or a \mobius band, respectively, with $n$ exceptional fibers. 

$M_K, M_L$ is the collection of the exteriors of 2-bridge knots and links, respectively.

$L_K$ is the collection of the exteriors of $(1,1)$-knots in lens spaces each of which admits a complete hyperbolic structure or admits a Seifert fibration whose regular fiber is not a meridian loop.

Let us examine the incompressible tori embedded in $M$ for each of the cases in the above theorem. 

\begin{lemma}\label{lem:onlyverticaltori}
Let $M$ be as in Theorem \ref{theo:kobayashi}. 
    If $T$ is an incompressible torus in $M$, then $T$ is either isotopic to a torus of the JSJ-decomposition, or to a vertical torus contained in a  Seifert fibered space in which $M$ is decomposed.
\end{lemma} 
\begin{proof}
    Note that the only tori, in each of the Seifert fibered spaces in which $M$ is decomposed, are those which are vertical surfaces. On the other hand, $M$ cannot contain  tori that come from two annuli, lying in different Seifert fibered spaces, glued along their boundary. If that were the case, the Seifert fibration of one of the pieces can be extended to the other one, possibly by changing the fibration, and giving a Seifert fibration for $M$, which is not possible. 
\end{proof}

\begin{cor} Let $M$ be as in Theorem \ref{theo:kobayashi}, then:

\begin{itemize}
    \item Type (i): contains exactly one incompressible separating torus.
    \item Type (ii):
    \begin{itemize}
        \item If $M_1 \in M\ddot{o}(n)$ with $n=0,1$, then there is only one separating incompressible torus. 
        \item If $M_1 \in M\ddot{o}(2)$ then there are three separating tori, two of them contained in $M_1$ and intersecting each other in four essential curves (see Figure \ref{Mo(2)}).
        \begin{figure}
            \centering
            \includegraphics[width=0.3\linewidth]{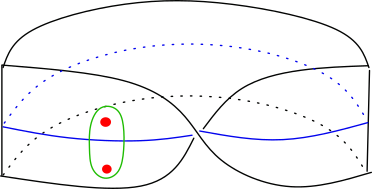}
            \caption{$M\ddot{o}(2)$. The green and blue circles lift to two incompressible separating tori, which intersect in four essential curves.}
            \label{Mo(2)}
        \end{figure}
    \end{itemize}
    \item Type (iii): 
    \begin{itemize}
        \item If $M_1 \in D(2)$ then there is only one incompressible separating torus.
        \item If $M_1\in D(3)$ then there are four incompressible separating tori. Three of them contained in $M_1$ and intersecting two by two in two essential curves (see Figure \ref{D(3)}).
        \begin{figure}
            \centering
            \includegraphics[width=0.3\linewidth]{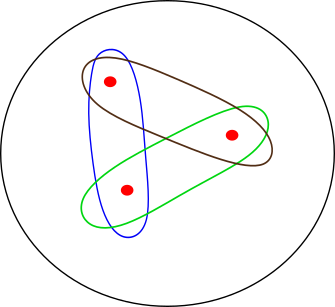}
            \caption{$D(3)$. The blue, green and brown circles lift to three incompressible, separating tori which intersect by pairs in two essential curves.}
            \label{D(3)}
        \end{figure}
    \end{itemize}
    \item Type (iv): There are two separating incompressible tori.
    \item Type (v):
    \begin{itemize}
        \item If $M_1\in A(n)$ with $n=0,1$, then there are two incompressible non-separating tori.
        \item If $M_1 \in A(2)$ then there are  two non-separating incompressible tori and $M_1$ contains two separating tori that intersect in four essential curves (see Figure \ref{A(2)}).
        \begin{figure}
            \centering
            \includegraphics[width=0.3\linewidth]{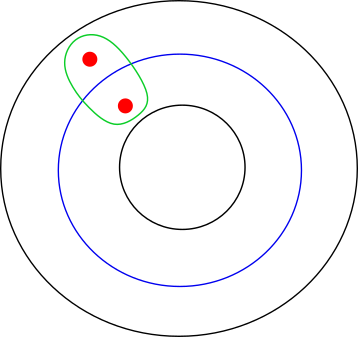}
            \caption{$A(2)$. The blue and green circles lift to two separating incompressible tori, which intersect in four essential curves.}
            \label{A(2)}
        \end{figure}
    \end{itemize}
\end{itemize}
\end{cor}
\begin{proof}
    By Lemma \ref{lem:onlyverticaltori} the possible incompressible tori in $M$ arises by lifting of essential curves that contain at least two singular points in the base space. Hence, by inspecting those curves we obtain the results. 
\end{proof}

\section{Bounds for the circular width of hyperbolic genus one knots} \label{bounds}

Let $K\subset \Sph^3$ be a non-fibered hyperbolic genus one knot, and $\Tmax \subset \calT$ a maximal simplicial collection. Suppose that the complementary regions of $\Tmax$ are handlebody pairs. Observe that a circular decomposition for $E(K)$ with circular width equal to $(3)$ can be achieved if there is a genus two Seifert surface $S$ for $K$ such that $(E(K), T, S)$ is a circular thin decomposition with $T\in \calT$. Hence $(3)$ is a lower bound for the circular width $cw(E(K))$.

In the case where $\vert\calT\vert=1$, the knot is a free genus one knot, and hence is almost fibered with circular width equal to $(3)$ or $(5)$, see \cite{FreegenusoneMNR}. 

We will assume that $\Tmax$ contains at least two tori.

\subsection{Lower bound for $cw(E(K)$}
 We will show that for many of these knots, the circular width is not equal to $(3)$.

 Let $K$ be a genus one hyperbolic knot, and let $\overline M$ be the corresponding manifold obtained from Construction \ref{construction}. If $cw(E(K))=(3)$, then by Corollary \ref{coro:theclosedmanifoldforextKisnotSF} $\overline M$ is not a Seifert fibered manifold. Moreover, it has a non-trivial torus decomposition and has a Heegaard splitting of genus two. 
Then Theorem \ref{theo:kobayashi} can be applied to $\overline M$. 

\begin{theorem}\label{dimMS(K)4o5}
    Let $K$ be a genus one hyperbolic knot in $\Sph^3$, $\Tmax \subset \calT$ a maximal simplicial collection, such that the complementary regions of $\Tmax$ are handlebody pairs. Suppose that $d= dim MS(K)= 4$ or $3$, then $E(K)$ does not have circular width equal to $(3)$.
\end{theorem}
\begin{proof}
   
     If $d=4$, then by Theorem \ref{theo:luisvaldez} $MS(K)$ is a $4$-simplex and $\calT= \Tmax$.  The corresponding collection of closed tori $\overline \Tmax$ consists of four disjoint, incompressible, separating tori in $\overline M$. But none of the Type(i)-(iv) of Theorem \ref{theo:kobayashi} contains those kinds of tori.
    
    If $d=3$, then by Theorem \ref{theo:luisvaldez} $MS(K)$ is 
    either a $3$-simplex, or $MS(K)$ consists of two $3$-simplices identified along a $2$-face.  In the first case, $\vert \calT \vert =4$, i.e., $\calT= \Tmax$ and the corresponding collection of tori $\overline \Tmax$ contains three disjoint incompressible, separating tori in $\overline M$, and none of the Types (i)-(v) of Theorem \ref{theo:kobayashi} contains those kinds of tori. If  $MS(K)$ consists of two $3$-simplices identified along a $2$-face, let $T_1, T_5 \in \calT$ be the two non-disjoint tori. Suppose that the circular thin decomposition for $E(K)$ has $T_1$ or $T_5$ as the thin level; therefore $\overline \calT$ consists of three disjoint, incompressible, separating tori in $\overline M$. If the thin level for the circular decomposition is not $T_1$ nor $T_5$, then $\overline \calT$ contains three disjoint, incompressible, separating tori and a fourth torus intersecting one of the former. In any of the cases, the manifolds of Types (i)-(v) cannot have such kinds of tori.
    
    Therefore, the circular width of  $E(K)$ is not equal to $(3)$.
\end{proof}

\begin{remark}
    For an infinite family of hyperbolic genus one knots with $MS(K)$ a $3$-simplex, the above result was proven in \cite{MarioAraEnrique}. Their proof is by contradiction, exploring the intersection of a thick surface of genus two with the family of Seifert tori, and the properties of the described family.  Moreover, the authors proved that the circular width, for that family of knots, is $(3,3)$. We include Figure \ref{family4toros} for reference.
\begin{figure}
    \includegraphics[scale=.72]{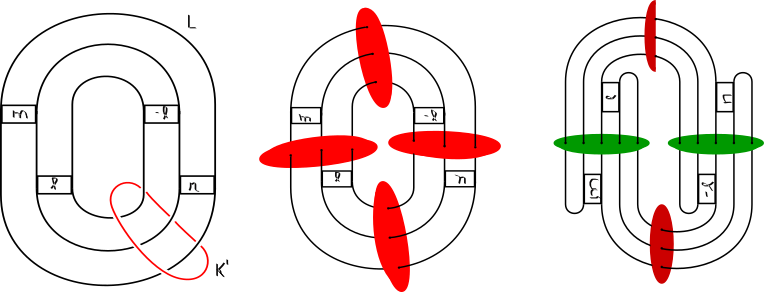}
   \caption{A family of knots with four non-isotopic Seifert surfaces [Cf.\cite{MarioAraEnrique}]}
    \label{family4toros}
\end{figure}
\end{remark}

\begin{theorem}
    Let $K$ be a genus one hyperbolic knot in $\Sph^3$, $\Tmax \subset \calT$ a maximal simplicial collection, such that the complementary regions of $\Tmax$ are handlebody pairs. Suppose that $d= dim MS(K)= 2$, then $E(K)$ does not have circular width equal to $(3)$ unless the manifold $\overline M$ obtained by Construction \ref{construction} is of Type (iv) (see Theorem \ref{theo:kobayashi}).
\end{theorem}
\begin{proof}

If $d=2$, then by Theorem \ref{theo:luisvaldez} $MS(K)$ is either a  $2$-simplex or $MS(K)$ consists of two $2$-simplices identified along one 1-face.

If $MS(K)$ is a $2$-simplex, then $\calT$ contains three disjoint Seifert tori, and the collection $\overline \calT$ in $\overline M$ consists of two disjoint, incompressible and separating tori. The cases Type (i),(ii),(iii) and (v) of Theorem \ref{theo:kobayashi} do not have these types of tori. Thus, $E(K)$ does not have circular width $(3)$ unless $\overline M$ is of Type (iv). 

If $MS(K)$ consists of two $2$-simplices identified along one 1-face, then $\calT$ contains four Seifert tori, two of which intersect each other in two essential curves, denote these tori by $T_1$ and $T_4$. Suppose that either $T_1$ or $T_4$ is the thin surface of the circular thin position for $E(K)$, then the collection $\overline \calT$ in $\overline M$ consists of two disjoint, incompressible and separating tori, and as before, $\overline M$ cannot be of Type (i),(ii),(iii) and (v) of Theorem \ref{theo:kobayashi}. Hence $E(K)$ does not have circular width $(3)$ unless $\overline M$ is of Type (iv).

On the other hand, if the thin surface of the circular thin position for $E(K)$ is different from $T_1$ and $T_4$, then $\overline \calT$ in $\overline M$ consists of three incompressible separating tori, two of them intersecting in two essential curves. The cases Type (i)-(v) of Theorem \ref{theo:kobayashi} do not have these types of tori. Thus, $E(K)$ does not have circular width $(3)$. 

\end{proof}

\begin{theorem}
    Let $K$ be a genus one hyperbolic knot in $\Sph^3$, $\Tmax \subset \calT$ a maximal simplicial collection, such that the complementary regions of $\Tmax$ are handlebody pairs. Suppose that $d= dim MS(K)= 1$, then $E(K)$ does not have circular width equal to $(3)$ unless the manifold $\overline M$ obtained by Construction \ref{construction} is of Type (i), Type (ii) with $M_1 \in M\ddot{o}(0,1)$, or Type (iii) with $M_1 \in D(2)$ (see Theorem \ref{theo:kobayashi}).
\end{theorem}
\begin{proof}
    
    If $d=1$, then by Theorem \ref{theo:luisvaldez} $MS(K)$ is either a $1$-simplex or $MS(K)$ consists of two $1$-simplices identified along a vertex.

    If $MS(K)$ is a $1$-simplex, then $\calT$ contains two disjoint Seifert tori, and the collection $\overline \calT$ in $\overline M$ consists of one incompressible and separating torus. The cases Type (ii) with $M_1 \in M\ddot{o}(2)$, Type (iii) with $M_1 \in D(3)$, Type (iv)-(v) of Theorem \ref{theo:kobayashi} do not have this type of torus. Hence $E(K)$ cannot have circular width $(3)$, unless $\overline M$ is of  Type (i), Type (ii) with $M_1 \in M\ddot{o}(0,1)$ or Type (iii) with $M_1 \in D(2)$

    If $MS(K)$ consists of two $1$-simplices identified along a vertex, then $\calT$ contains three Seifert tori two of which intersect each other in two essential curves, denote these tori by $T_1$ and $T_3$. Suppose that either $T_1$ or $T_3$ is the thin surface of the circular thin position for $E(K)$, then the collection $\overline \calT$ in $\overline M$ consists of one incompressible and separating torus, and as in the above case, $E(K)$ does not have circular width $(3)$  unless $\overline M$ is of  Type (i), Type (ii) with $M_1 \in M\ddot{o}(0,1)$ or Type (iii) with $M_1 \in D(2)$.
    
      If the thin surface of the circular thin position for $E(K)$ is different from $T_1$ and $T_3$, then $\overline \calT$ in $\overline M$ consists of two incompressible, separating tori that intersect each other in two essential curves. The cases of Type (i)-(v) of Theorem \ref{theo:kobayashi} do not have these types of tori. Thus, $E(K)$ does not have circular width $(3)$. 
    
\end{proof}

\subsection{Upper bound for $cw(E(K))$}

Now we give some upper bounds for the circular width of a genus one hyperbolic knot. In particular, we consider the case where all complementary regions are simple handlebody pairs.

\begin{lemma}\label{regionR} 
Let $K$ be a genus one hyperbolic knot in $\Sph^3$, $\Tmax \subset \calT$ a maximal simplicial collection. Let $R=R_{i,i+1}$ be the complementary region between the tori $T_i$ and $T_{i+1}$. Suppose that $(R,J)$ is a handlebody pair. Then $R$ has a decomposition $(A_i,B_i;S_i)$ of complexity $c(S_i)=5$.

    \end{lemma}

    \begin{proof}
     By hypothesis $R$ is a genus two handlebody and $J\subset R$ is a non-trivial separating curve. Let $H \subset \textrm{int} R$ a genus two handlebody such that the complement of $H$ in $R$ is a product $\bdry R \times I$. Let $\tau$ be a straight arc in $\bdry R \times I$ connecting $T_i$ and $H$. Note that $N(\tau)\cup H$ is equivalent to adding two $1$-handles to $T_i\times [0,\frac{1}{2}]$, along $T_i\times\{\frac{1}{2} \}$, since $H$ is a genus two handlebody. Then $A_i=T_i\times [0,\frac{1}{2}]\cup N(\tau)\cup H$ is a compression body with $\bdry_- A_i= T_i$ and $\bdry_+ A_i= T_i \# \bdry H$ which is a genus three Seifert surface for $K$. Let $\alpha_1$, $\alpha_2$ be two essential non-parallel simple closed curves in $T_i\times \{\frac{1}{2}\}$, which intersect precisely in $\bdry \tau$. Note that the complement in $T_i\times\{\frac{1}{2}\}$ of the union of $\alpha_1$ and $\alpha_2$ is an annulus which is a neighborhood of the curve $J$. Now, $\alpha_i\times [\frac{1}{2},1] -\textrm{int}N(\tau)$, is a disk whose boundary consists of one arc in $T_i\times\{\frac{1}{2}\}$, one arc in $\bdry H$ and two arcs in $\partial N(\tau)$.  A neighborhood of these disks can be considered as two 2-handles, $h_1, h_2$, attached to $A_i$ along $\bdry_{+}A_i$. The boundary of the resulting manifold is a torus parallel to $T_{i+1}$. In fact, $\bdry_+ A_i \times I \cup h_1 \cup h_2$ is a compression body $B_i$ with $\bdry_+B_i=\bdry_+ A_i$ and $\bdry_-B_i=T_{i+i}$. Note that the surface $S_i=T_i\# \partial H$ is a thick surface of genus $3$, with complexity $c(S_i)=1-\chi(\bar{S_i})=5$.
     \end{proof}

\begin{lemma}\label{regionsimple} 
Let $K$ be a genus one hyperbolic knot in $\Sph^3$, $\Tmax \subset \calT$ a maximal simplicial collection. Let $R=R_{i,i+1}$ the complementary region between the tori $T_i$ and $T_{i+1}$. Suppose that $(R,J)$ is a simple pair.  Then $R$ has a decomposition $(A_i,B_i;S_i)$ of complexity $c(S_i)=3$. Furthermore, if there are two adjacent simple pairs $R=R_{i,i+1}$ and $R'=R_{i+1,i+2}$ such that the corresponding power curves intersect in one point, then $R\cup R'$  has a decomposition $(A,B;S)$ of complexity $c(S)=3$. 
\end{lemma}

\begin{proof} 
Suppose now that $R$ is a simple pair. Let $V$ be the companion solid torus contained in $R$. Note that there is an annulus $C$ contained in $R$, such that one boundary component lies on $T_i$ and the other in $\bdry V$. Let $\tau$ be a straight arc contained in $C$ connecting $T_i$ and $\bdry V$. Then $N(\tau)\cup V$ can be considered as a 1-handle attached to $T_i\times [0, \frac{1}{2}]$, and $A_i= T_i\times [0,\frac{1}{2}] \cup N(\tau)\cup  V$ is a compression body, with $\bdry_{-}A_i = T_i$ and $\bdry_{+} A= T_i\# \bdry V$, which is a genus two Seifert surface for $K$. Let $D=C-\textrm{int} N(\tau)$, this is a disk whose boundary consists of an arc lying in $T_i$, two arcs in $\bdry N(\tau)$ and one arc in $\bdry V$. A regular neighborhood of this disk can be considered as a 2-handle $h$ attached to $\bdry_+A_i$.
In fact, $\bdry_+ A_i \times I \cup h$ is a compression body $B_i$ with $\bdry_+B_i=\bdry_+ A_i$ and $\bdry_-B_i=T_{i+i}$. Note that the surface $S_i=T_i\# \partial H$ is a thick surface of genus two, with complexity $c(S_i)=1-\chi(\bar{S_i})=3$.

Suppose now that $R=R_{i,i+1}$ and $R'=R_{i+1,i+2}$ are two consecutive simple pairs. Let $V$ and $V'$ be the companion solid tori contained in $R$ and $R'$, respectively. Suppose that $V$ and $V'$ are large enough such that $V\cap T_{i+1}$ and $V'\cap T_{i+1}$ each consists  of an annulus, and such that the annuli intersect in one disk. This is equivalent to saying that the corresponding power curves intersect in a point. Suppose also that $\bdry V$ contains two spanning annuli of $R_i$, and similarly for $\bdry V'$. Let $\tau$ be a straight arc connecting $T_i$ and $\bdry V'$, disjoint from $V$, that is, an arc contained in a spanning annulus of $R_i$, but pushed off $V$. Then $N(\tau)\cup V'$ can be considered as a 1-handle attached to $T_i$. Then $A= T_i\times [0,\frac{1}{2}] \cup N(\tau)\cup  V'$ is a compression body, with $\bdry_{-}A= T_i$ and $\bdry_{+} A_i= T_i\# \bdry V'$, which is a genus two surface for $K$. Similarly, let $\tau'$ be a straight arc connecting $T_{i+2}$ and $\bdry V$, disjoint from $V'$. Then $N(\tau')\cup V$ can be considered as a 1-handle attached to $T_{i+2}$, and $B= T_{i+2}\times [0,\frac{1}{2}] \cup N(\tau)\cup  V$ is a compression body, with $\bdry_{-}B = T_{i+2}$ and $\bdry_{+} B = T_{i+2}\# \bdry V'$, which is a genus two surface for $K$. Note that $\bdry_{+}A$ and $\bdry_{+}B$ are parallel Seifert surfaces for $K$, which gives a decomposition of $R$, with complexity $3$.
\end{proof}

Let $(H(p/q),J)$  be a simple pair and let $\gamma_0$ and $\gamma_1$ be curves in $T\times \{0 \}$ and $T\times \{ 1\}$, respectively,  which intersect a power curve $\omega$ of $H(p/q)$ in one point. Take a basis for $H_1(T)$ consisting of $[\omega]$ and $[\beta ]$, where $\beta$ is a curve that intersects $\omega$ in one point. Then $[\gamma_0] = [\beta] +n_0[\omega]$ and $[\gamma_1]=[\beta] + n_1[\omega]$. We define the \emph{$\omega$-linking number of $\gamma_0$ and $\gamma_1$},
as $\vert n_1-n_0 \vert $.

    \begin{lemma}\label{regionsimple3} 
    
    Let $K$ be a genus one hyperbolic knot in $\Sph^3$, $\Tmax \subset \calT$ a maximal simplicial collection. Suppose that there are three adjacent regions, say $R_1$, $R_2$, $R_3$, bounded by the tori $T_1$, $T_2$, $T_3$ and $T_4$. Suppose that $(R_i, J_i)$ is a simple pair for each $i=1,2,3$. Assume that the corresponding power curves in adjacent regions intersect in a point, and that the power curves of $R_1$ and $R_3$ have $\omega$-linking number $n$ in $R_2$.
    
    Then the following are equivalent:
    \begin{enumerate}
        \item $R=R_1\cup R_2 \cup R_3$ has a decomposition of complexity $3$.
        \item $R_2=(H(a/b),J)$  for some  $(a,b)$ such that $b-an=\pm 1$.
        \item $R_2[J]$ is the exterior of a torus link $(2,2p)$ and the power curves of $R_1$ and $R_3$ are meridian curves of $R_2[J]$. 
    \end{enumerate}
     
    \end{lemma}
   
\begin{proof} First we prove $(1)$ implies $(3)$. Apply the third step of  Construction \ref{construction} to $R$. The resulting manifold $\overline R$ is a graph manifold containing two incompressible tori, $\overline{ T_2}$ and $\overline{ T_3}$, which divide $\overline R$ into two Seifert fibered spaces over a disk with two exceptional fibers and a Seifert fibered space over an annulus, namely $\mathbb{A} (b/a)$. If $R$ has a decomposition of complexity $3$ then $\overline R$ will have Heegaard genus two, and case (iv) of  Theorem \ref{theo:kobayashi}  implies that the central piece, $\mathbb{A}(b/a)$ must be the exterior of a 2-bridge link, and that the power curves of $R_1$ and $R_3$ are identified to meridians of the link. Now, the only 2-bridge links whose exterior is a Seifert fibered space are the $(2,2p)$-torus links, whose exterior has a fibration of type $\mathbb{A}(1/p)$. This shows that if $R$ has a decomposition of complexity $3$ then $R_2[J]$ is the exterior of a torus link $(2,2p)$ and the powers curves of $R_1$ and $R_3$ are meridians of $R_2[J]$. This proves $(1)$ implies $(3)$.

Now we show that $(2)$ implies $(1)$. Suppose that $R_2=(H(a/b),J)$ for some $(a,b)$ such that $b-an=\pm 1$. We show that $R$ has a decomposition of complexity $3$. To do that, we construct an explicit decomposition.

Consider $R_1=(T_1\times [0,1])\cup V_1$, where $T_1=T_1 \times \{0\}$, and $V_1$ is a solid torus such that $(T_1\times [0,1])\cap V_1$ is an annulus $A_1'$, contained in $T_1\times \{1\}$ and $\bdry V_1$, which is essential in both $T_1\times [0,1]$ and $V_1$. Let $A_1$ be the complementary annulus of $A_1'$ in $\bdry V_1$. Note that $T_2= (T_1\times \{1\} - A_1') \cup A_1$. Similarly, let $R_3=(T_4\times [0,1])\cup V_3$, where $T_4=T_4 \times \{1\}$, and $V_3$ is a solid torus such that $(T_4\times [0,1])\cap V_3$ is an annulus $A_3'$, contained in $T_4\times \{0\}$ and $\bdry V_3$, which is essential in both $T_4\times \{0\}$ and $V_3$. Let $A_3$ be the complementary annulus of $A_3'$ in $\bdry V_3$. Note that $T_3= (T_4\times \{0\} - A_3') \cup A_3$. Consider now $T_2\times [0,1]$, and let $A_2$ be an essential annulus in $T_2$. Let $V=A_2\times [0,1] \subset T_2 \times [0,1]$. Consider  $R_2 = (T_2\times [0,1] - \textrm{int} V) \cup V_2$, where $V_2$ is a solid torus such that $\bdry V_2$ is identified with $\bdry V$. Note that $T_2=T_2\times \{0\}$ and $T_3=T_2\times \{1\}$.

Note that the cores of $A_1$ and $A_3$ represent the power curves of $R_1$ and $R_3$, respectively. We assume that the core of $A_1$ and $A_2\times \{0\}$ (both contained in $T_2$) intersect in one point, and that the core of $A_3$ and the core of $A_2\times \{1\}$ (both contained in $T_3$) intersect in one point. Taking $A_1$ and $A_3$ very thin and $A_2$ to be fat, we can assume that $A_1$ is the union of two rectangles, one in $T_2-(A_2\times \{0\})$ and the other inside $A_2\times \{0\}$, and that $A_3$ is the union of two rectangles, one in $T_3-(A_2\times \{ 1 \})$ and the other in $A_2 \times \{1\}$. Isotopying part of $A_3$ inside $A_2\times \{1\}$, we can suppose that the part of $A_1$ lying in $T_2-(A_2\times \{ 0 \})$ looks parallel to the part of $A_3$ lying in $T_3-(A_2\times \{ 1 \}$. See Figure \ref{fig:9}.
But the part of $A_1$ lying inside $A_2\times \{ 0 \}$ may not look parallel to the part of $A_3$ lying inside $A_2\times \{ 1 \}$. In fact, the difference in twists is exactly the number $n$ of intersections between the core curves of $A_1$ and $A_3$.

\begin{figure}[h]
\centering
\includegraphics[scale=0.5]{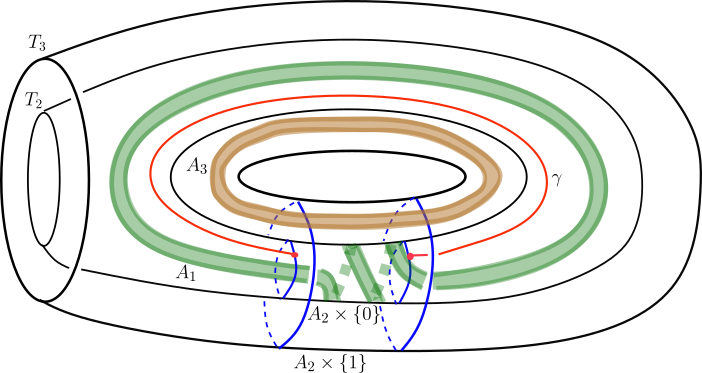}
\caption{The region $R_2$.}
\label{fig:9}
\end{figure}

\begin{figure}[h]
\centering
\includegraphics[scale=0.5]{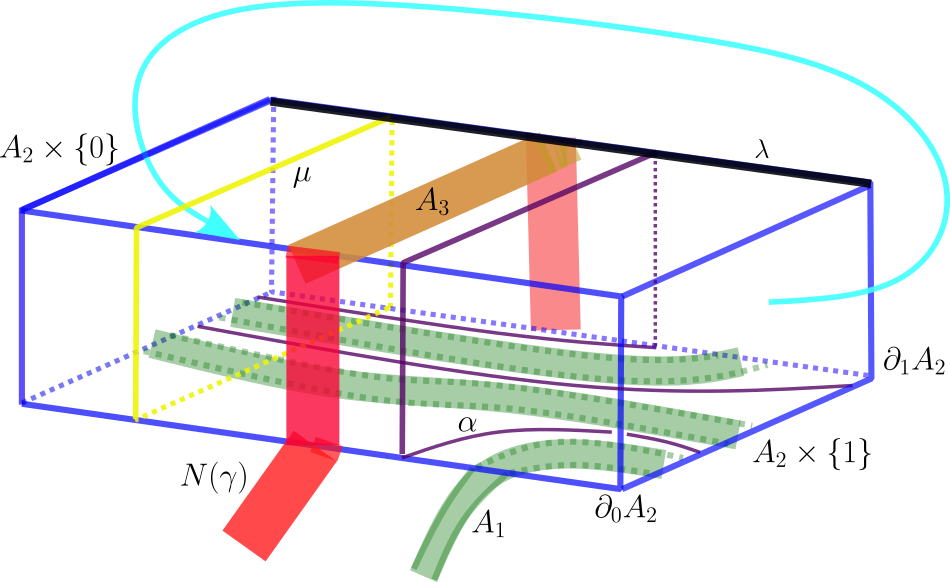}
\caption{The region around the solid torus $V_2$.}
\label{fig:10}
\end{figure}

Let $\bdry_0 A_2$ and $\bdry_1 A_2$ denote the boundary components of $A_2$. $\bdry_0 A_2\times [0,1]$ and $\bdry_1 A_2\times [0,1]$ are two spanning annuli in $R_2$. Let $\gamma$ be an arc in $T_2-\textrm{int} A_2$, with endpoints in $\bdry_0 A_2$ and $\bdry_1 A_2$, disjoint from $A_1$, and running parallel to a component of $\bdry A_1$. Let $N(\gamma)$ be a neighborhood of $\gamma$ in $T_2-\textrm{int} A_2$, disjoint from $A_1$. See Figure \ref{fig:10}. Let $N(\gamma)\times [0,1]$ be a 3-ball such that $N(\gamma)\times \{1\} = A_3-(\textrm{int} A_2 \times \{ 1 \})$.  $N(\gamma)\times [0,1]$ intersects each of these annuli in a disk. Consider $H_1=T_1\times [0,1] \cup (N(\gamma)\times [0,1]) \cup V_3$, which is clearly a compression body, for it is the disk sum of a product and a solid torus. Note that $\bdry_-H_1=T_1\times \{ 0\}$ and $\bdry_+H_1= T_1\times \{ 1\}\sharp \bdry V_3$ which is a genus two surface $S$ with boundary parallel to $J$.

The surface $S$ splits $R$ into two pieces, one being $H_1$. It remains to see that the complementary part $H_2$ is also a compression body. To see this, first consider $T_4\times [0,1]$. Then consider $A''=(T_3-\textrm{int} A_3)-\textrm{int} A_2$. This is just an annulus. Now, let $H'=T_4\times [0,1]\cup (A''\times [0,1])$, where  $T_4\times [0,1]$ is contained in $R_3$, while $(A''\times [0,1])$ is contained in $R_2$. Note that $H'$ is homeomorphic to $T_4\times [0,1]$. Note that  $H'\cap V_2$  consists of three rectangles: one  contained in $A_2\times \{1\}$, another in $\bdry_0 A_2\times [0,1]$ and one more in $\bdry_1 A_2\times [0,1]$. The union of the three rectangles is a disk $C$. See Figure \ref{fig:11}.

\begin{figure}[h]
    \centering
\includegraphics[scale=0.5]{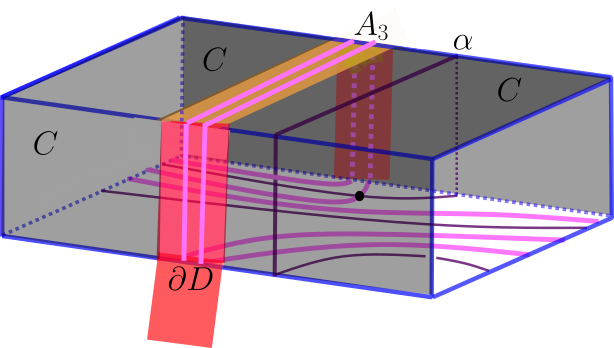}
    \caption{The region around the solid torus $V_2$.}
    \label{fig:11}
\end{figure}

 Then $H'\cup V_2$ is a compression body, as it is the disk sum of a product and a solid torus. 

Observe that $H_2$, the complementary part of $H_1$, is $H_2=H'\cup V_2 \cup V_1$, and $(H'\cup V_2)\cap V_1$ consists of an annulus, that is, it is precisely $A_1$. But $A_1$ intersects $H'$ and $V_2$, each one on a disk. As $A_1$ goes around $V_1$ at least twice longitudinally, $H_2$ will be a compression body if only if $A_1$ intersects a meridian of $V_2$ in a single arc.  

Let $\alpha$ be a curve in $\bdry V_2$ that consists of an arc in $A_2\times  \{ 0 \}$, which runs parallel to $A_1$, an arc in $A_2 \times \{ 1 \}$, which runs parallel to $A_3$, and two straight arcs in $\bdry_0 A_2 \times [0,1]$ and $\bdry_1 A_2 \times [0,1]$. Taking a natural meridian $\mu$ and longitude $\lambda$ for the solid torus $V$, it is clear that the curve $\alpha$ has slope $(1,n)$ (that is, one meridian and $n$ longitudes). Let $D$ be a meridian disk of $V_2$. We can assume that $\partial D$ lies in $\partial V_2 - C$, and then it intersects the part of $A_1$ lying in $\partial V_2$ in a single arc.  This is the same as saying that it intersects the curve $\alpha$ in one point, and then $\bdry D$ has slope $(a,b)$, where $b-an=\pm 1$. See Figure \ref{fig:11}. This shows that $H_2$ is a compression body if $R_2= H(a/b)$, where $b-an=\pm 1$. Therefore, $R$ has a decomposition $(H_1, H_2; S)$ with complexity $3$ if $R_2= (H(a/b), J)$, where $b-an=\pm 1$.

It remains to show that $(3)$ implies $(2)$: that $R_2=(H(a/b),J)$ for some $(a,b)$ such that $b-an= \pm 1$ if $R_2[J]$ is the exterior of a torus link $(2,2p)$ and the power curves of $R_1$ and $R_3$ are meridians of $R_2[J]$.

Embed  $T\times[0,1]$ in $S^3$ in a standard way, let $L_2$ be a meridian of $T$ lying in the level $\{1/2\}$. By performing $a/b$-Dehn surgery on $L_2$,  $T\times[0,1]$ becomes $R_2[J]$. Let $L_1$ and $L_3$ be the cocores of the complementary solid tori of the embedding of $T\times[0,1]$ in $S^3$, where $L_1$ is the black curve shown in Figure \ref{fig:12}$(i)$, while $L_3$ is the green curve and the blue curve represents $L_2$. Then, up to homeomorphism, $R_2[J]$ is homeomorphic to the exterior of the link $L_1 \cup L_2 \cup L_3$, after performing $a/b$- Dehn surgery on $L_2$. The slopes $n$ and $1/0$ on $L_1$ and $L_3$ represent the slopes of the power curves in $R_1$ and $R_3$, respectively, but we are not performing surgery on these components. Now, perform $(-n-1)$-twists around $L_2$; the result is the link shown in Figure \ref{fig:12}$(ii)$. 
Next, do a twist around $L_1$ and apply an isotopy, to obtain the link shown in Figure \ref{fig:12}$(iii)$. Note that $(b-an)/(bn-an-a)$-Dehn surgery on $L_2$ produces the exterior of a link in $S^3$ (with the slopes $1/0$ as meridians) if and only if $b-an=\pm 1$. Suppose that we have that condition. We get the link shown in Figure \ref{fig:12}$(iv)$. Now, by performing $1/(1-a)$-twists around $L_2$, we get the link shown in Figure \ref{fig:12}$(v)$, which is the exterior of the $(2,2a)$-torus link.

\begin{figure}[h]
    \centering
\includegraphics[scale=0.5]{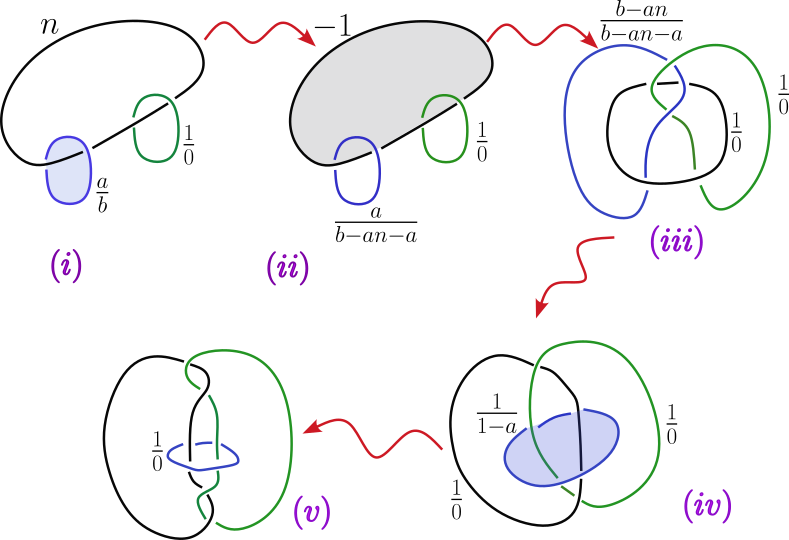}
    \caption{$H(a/(na\pm 1))[J]$ is the exterior of the $(2,2a)$-torus link}
    \label{fig:12}
\end{figure}

\end{proof}

\begin{theorem}\label{width5}
    Let $K$ be a genus one hyperbolic knot in $\Sph^3$, $\Tmax \subset \calT$ a maximal simplicial collection, such that the complementary regions of $\Tmax$ are handlebody pairs. Let $d= dim MS(K)$.
    
    Then $E(K)$ has a circular decomposition $(E(K), \calF, \calS)$ of circular width $cw(E(K), \calF, \calS) =(5,5,...,5)$, which is a $(d+1)$-tuple.

    If there are $ 1 \leq l \leq d+1$ complementary regions which are simple pairs, then  $E(K)$ has a circular decomposition $(E(K), \calF, \calS)$ of circular width a $(d+1)$-tuple with $'5'$ in the first $d+1-l$ entries and $'3'$ in the last $l$ entries.
    
     If there are two adjacent simple pairs and the corresponding power curves intersect in one point, then $E(K)$ has a circular decomposition with circular width a $d$-tuple with $'5'$ in the first $d_1$ entries, $'3'$ in the following $d_2$ entries and $'3'$ in the last entry, where  $d_1+d_2=d-1$.
\end{theorem}

\begin{proof}

By hypothesis the complementary regions $R_{i,i+1}$ of $\Tmax$ are handlebody pairs, for $i=1,...,d+1$. Lemma \ref{regionR} implies that each region has a decomposition $(A_i,B_i; S_i)$ with $c(S_i)=5$. Let $\calF$  be equal to $\Tmax$ and $\calS= \{S_i\}_{i=1}^{d+1}$, then  $(E(K), \calF, \calS)$ is a circular decomposition with circular width $(5,5,...,5)$ a $(d+1)$-tuple.

If a region $R_{i,i+1}$ is a simple pair, then Lemma \ref{regionsimple} implies that such a pair has a decomposition with complexity $3$. Now suppose that $l$ regions are simple pairs, then $E(K)$ has a circular decomposition with $l$ thick surfaces of complexity $3$ and $d+1-l$ thick surfaces of complexity $5$. Hence the circular width of the decomposition is a $(d+1)$-tuple with '$5$' in the first $d+1-l$ entries and '$3$' in the last $l$ entries.

If there are two adjacent simple pairs $R=R_{i,i+1}$ and $R'=R_{i+1,i+2}$ such that the corresponding power curves intersect in one point, then $R\cup R'$  has a decomposition $(A,B;S)$ of complexity $c(S)=3$, by Lemma \ref{regionsimple}. This implies that $E(K)$ has a circular decomposition with one thick surface of complexity $3$ coming from $R\cup R'$,  $d_1$ thick surfaces of complexity $5$ and $d_2$ thick surfaces of complexity $3$, with $d_1+d_2=d-1$.

\end{proof}

\begin{cor}\label{cor:cotasparaelancho}
     Let $K$ be a genus one hyperbolic knot in $\Sph^3$, $\Tmax \subset \calT$ a maximal simplicial collection, such that the complementary regions of $\Tmax$ are handlebody pairs, and $d= dim MS(K)$. 
     
     Then the circular width of $E(K)$ satisfies: 
     \begin{center}
         $(3)\leq cw(E(K))\leq (5,5,...,5)$
     \end{center}
     If all the complementary regions $\Tmax$ are simple pairs, then the circular width of $E(K)$ satisfies: 
       \begin{center}
         $(3)\leq cw(E(K))\leq (3,3,...,3)$
     \end{center}

     The right hand side of the inequalities is a $(d+1)$-tuple.
\end{cor}
\begin{proof} This is a direct consequence of Theorem \ref{width5}.
\end{proof}
\begin{cor}
     Let $K$ be a genus one hyperbolic knot in $\Sph^3$, $\Tmax \subset \calT$ a maximal simplicial collection. If all complementary regions of $\Tmax$ are simple handlebody pairs and $(3)<cw(E(K))$, then $K$ is non-almost fibered. 
\end{cor}
\begin{proof}
    If $K$ were almost fibered it would have circular width a $1$-tuple, which cannot be $(3)$ by hypothesis, so it must be $(n)$ with $n>3$. But also $(n)\leq (3,3,...,3)$, by Corollary \ref{cor:cotasparaelancho}, implying that $n\leq 3$. This contradiction shows that $K$ is non-almost fibered. 
\end{proof}

\subsection{Genus one knots with simple pairs and power curves intersecting in one point}

In the following, we prove that genus one knots with a disjoint collection of Seifert tori such that the complementary regions are simple pairs with the power curves of every two adjacent regions intersecting in a single point are hyperbolic knots. Moreover, we describe their Kakimizu complex and their circular width is computed.

\begin{theorem} \label{theo:MaximalKakimizuSimpleCircularwidth}
Let $K$ be a genus one knot that contains a collection of disjoint non-isotopic Seifert tori $T_1, T_2, \dots T_n$, which decomposes $E(K)$ into simple pairs. Suppose that for every two adjacent pairs, the corresponding power curves intersect in one point. Then

\begin{enumerate}
   
    \item $K$ is a hyperbolic knot or $K$ is the trefoil and $n=1$.
    \item There is no other Seifert torus for $K$.  
    \item The Kakimizu complex of $K$, $MS(K)$, is a $(n-1)$-simplex with $n\leq 5$.
    \item If $n=1$, then $K$ is fibered or $cw(E(K))=(3)$. 
    \item If $n=2$, then $cw(E(K))=(3)$
    \item If $n=3$, $cw(E(K))=(3)$ or $(3,3)$. 
    \item If $n=4$, $cw(E(K))=(3,3)$.
    \item If $n=5$, $cw(E(K))=(3,3)$ or $(3,3,3)$.
\end{enumerate}
\end{theorem}
\begin{proof}
$K$ is hyperbolic or $K$ is the trefoil and $n=1$ since the hypothesis satisfies [Cf.\cite{LVgenusonehypknotatmostfive} Lemma 8.1 (a)].
The proof of Lemma 4.8 in \cite{MarioAraEnrique} can be replicated in this case, hence $K$ cannot contain another genus one Seifert surface non-isotopic to any $T_i$. Since $K$ has exactly $n$ disjoint non-isotopic Seifert tori by $(2)$, then $MS(K)$ is an $(n-1)$-simplex, and since $K$ is hyperbolic Theorem \ref{theo:ValdezAtMost5tori} implies that $n\leq 5$. If $K$ is non-fibered and $n=1,2$, by Lemma \ref{regionsimple} $c(E(K))=(3)$.  In the case where $n=3$, if the conditions of Lemma \ref{regionsimple3} are satisfied then $cw(E(K))=(3)$, otherwise by Lemma \ref{regionsimple} $cw(E(K))=(3,3)$. If $n=4,5$, Theorem \ref{dimMS(K)4o5} implies that  $cw(E(K))>(3)$. By Lemma \ref{regionsimple} for $n=4$ $cw(E(K))=(3,3)$. Finally, for $n=5$, if the conditions of Lemma \ref{regionsimple3} are satisfied then $cw(E(K))=(3,3)$, otherwise, by Lemma \ref{regionsimple} $cw(E(K))=(3,3,3)$.
    
\end{proof}

\section{Examples}
\label{examples}

In this section, we give several families of examples satisfying Theorem \ref{theo:MaximalKakimizuSimpleCircularwidth}. The examples are constructed via double branched covers. In fact, a simple pair arises as the  double branch cover of a certain tangle. We give examples of knots having 3 or 5 disjoint non-parallel Seifert Tori. Example with 4 Seifert tori are given in \cite{MarioAraEnrique}, see Figure \ref{family4toros}.

Let $H'=D\times I$, where $D$ is a disk, let $J'=\bdry D \times \{ 1/2 \}$. Consider $\{x_1^i,x_2^i,x_3^i \}$ a collection of three points in the interior of $D\times \{ i \}$, for $i=0,1$. Let $\tau$ be a collection of three arcs properly embedded in $D\times I$, with endpoints the $x_j^i$'s. We denote this tangle by $(H',J',\tau)$.

\begin{figure}[h]
    \centering
\includegraphics[scale=0.5]{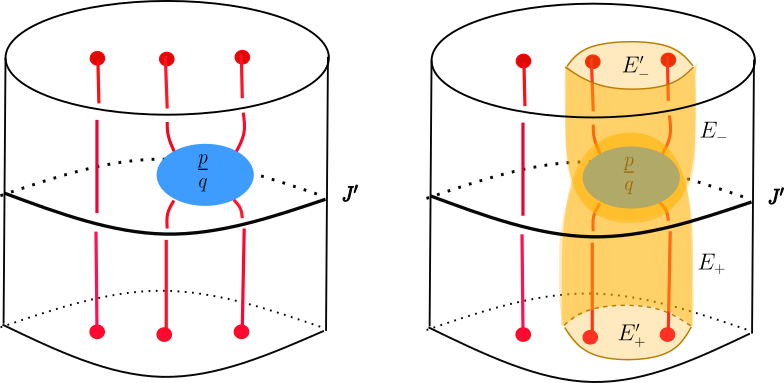}
     \caption{A simple-tangle}
    \label{simple-tangle}
\end{figure}

We say that a tangle $(H',J',\tau)$ is \emph{simple} if it is homeomorphic, up to a homeomorphism fixing $J'$, to a tangle where $\tau$ consists of a product arc and two arcs that form a rational tangle, see Figure \ref{simple-tangle}. 

\begin{lemma}\label{branchedcover} A simple pair $(H,J)$ double branches cover a simple tangle $(H',J',\tau)$. Conversely, the double cover of a simple tangle $(H',J',\tau)$ branched along $\tau$ is a handlebody simple pair $(H,J)$, where the lifting of $J'$ is the separating curve $J$.

\end{lemma}

\begin{figure}[h]
    \centering
\includegraphics[scale=0.6]{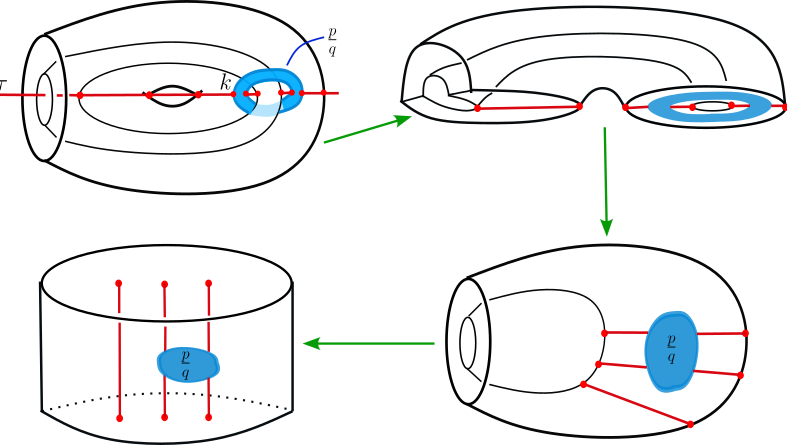}
     \caption{A simple pair covers a simple tangle}
    \label{cover-simple-tangle}
\end{figure}

\begin{proof} Let $T$ be a once punctured torus, and consider $T\times I$. Let $J= \bdry T \times \{ 1/2 \}$, $\gamma$ be a non-separating  curve in $T$ and let $k=\gamma \times \{ 1/2 \}$. Let $H(p/q)$ be the 3-manifold obtained by performing $p/q$-Dehn surgery on $k$. Note that $(H(p/q), J)$ is a simple pair, and in fact any simple pair is obtained in this way. Up to homeomorphism, we can assume that $\gamma$ is a meridian of $T$, as shown in Figure \ref{cover-simple-tangle}. Note that there is an involution $\tau'$ on $H(p/q)$, whose set of fixed points consists of three arcs. By taking the quotient space of this involution, we get a simple tangle $(H',J',\tau)$ as shown in Figure \ref{cover-simple-tangle}, where $J'$ is the image of $J$, and the arcs $\tau$ are the image of the arcs fixed by the involution $\tau'$.

Conversely, given $(H',J',\tau)$, we can construct its double branched cover from Figure \ref{cover-simple-tangle}, just going in the opposite direction.

\end{proof}

Let $(H_i',J',\tau)$ be a simple tangle, note that there are disks $E_{+}'$, $E_{+}$,  $E_{-}'$ and $E_{-}$, which satisfy the following:

\begin{itemize}
        \item $E_{+}'\subset D\times \{0\}$ and $E_{+}' \cap \tau=$ two points.
        \item $E_{+}$ is properly embedded in $H'_i$ and $E_{+} \cap \tau=$ two points.
        \item $\bdry E'_{+} =\bdry E_{+}$ 
        \item $E_{+}'\cup E_{+}$ bounds a 3-ball $B$ in $H_i'$ that determines the rational tangle $R(p/q)$.
     \item $E_{-}'\subset D\times \{ 1\}$ and $E_{-}' \cap \tau=$ two points.
        \item $E_{-}$ is property embedded in $H_i'$ and $E_{-} \cap \tau=$ two points.
        \item $\bdry E'_{-} =\bdry E_{-}$ 
        \item $E_{-}'\cup E_{-}$ bounds a 3-ball $B'$ in $H_i'$ that also determines the rational tangle $R(p/q))$.
    \end{itemize}
        
    \begin{remark}\label{discosenH} Disks such as $E_{\pm}'$ are illustrated in Figure \ref{simple-tangle}, these are the shaded disks in the boundary of the $3$-ball shown there. They are called companion disks of $H'$. The disks $E_{\pm}$ are properly embedded in $H'$, and their boundary coincide with that of $E_{\pm}'$. The union $E_{\pm}\cup E_{\pm}'$ is a sphere that bounds a $3$-ball which encloses the tangle $R(p/q)$, see Figure \ref{simple-tangle}. The liftings $C_{\pm}$ and $C_{\pm}'$ of $E_{\pm}$ and $E_{\pm}'$ are annuli. 
    The core of $C_{\pm}'$ is a power curve for $H_i$.

Suppose that two simple tangles $(H_1',J_1',\tau_1)$ and $(H_2', J_2',\tau_2)$ are identified along $D\times \{1\}$ in  $H_1'$ and $D\times \{0\}$ in $H_2'$.  Assume that there are companion disks $E_{i,\pm}\subset H_i'$, for $i=1,2$, such that $\bdry E_{1,-} \cap \bdry E_{2,+}$ consists of two points, then the lift of $E_{1,-}$ and $ E_{2,+}$ in the double branched covering of $(H_1',J_1',\tau_1)\cup (H_2', J_2',\tau_2) $ consists of two annuli whose cores intersect in a single point. 
\end{remark}

\subsection{A family of knots with three disjoint non-isotopic Seifert surfaces of genus one and circular width (3)}

\begin{figure}[h]
    \centering
    \def\svgwidth{1\linewidth}
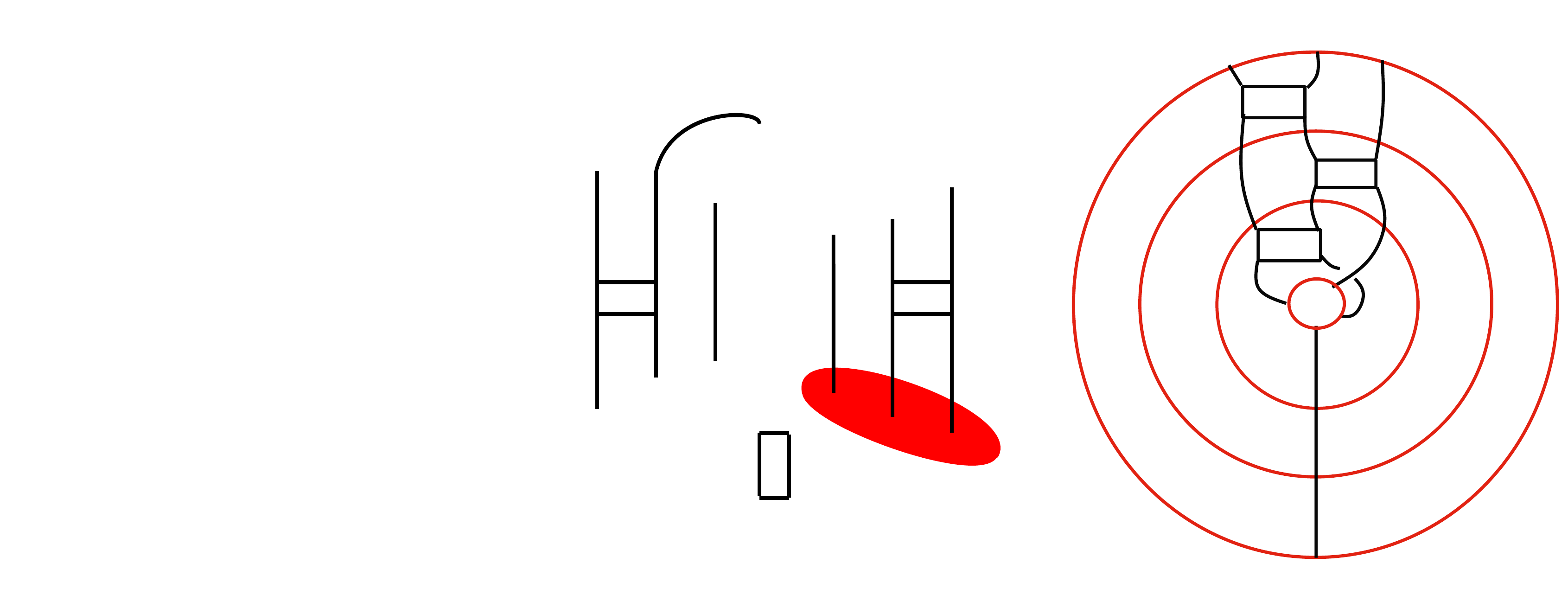
    \caption{A family of knots with three non-isotopic Seifert surfaces and circular width $(3)$.}
    \label{family3torosancho3}
\end{figure}

        Let $K'\cup L$  be the link shown on the left hand side of Figure \ref{family3torosancho3}. The boxes represent a sequence of  $l$ or $m$ horizontal crossings. Let us assume that $\vert l\vert, \vert m\vert \geq 2$. Both $K'$ and $L$ are trivial knots. This link is just a sample; more complicated examples can be found by choosing appropriate rational tangles instead of the integral tangles shown in the boxes.

\begin{figure}
    \centering
    \includegraphics[width=0.5\linewidth]{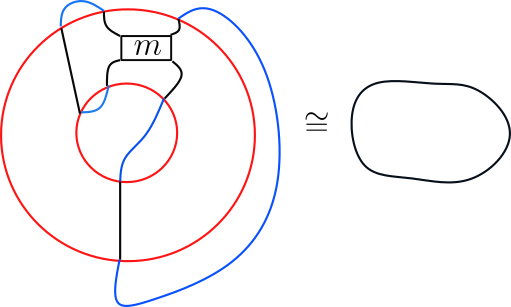}
    \caption{A filling that produces the trivial knot}
    \label{fillingtrivial}
\end{figure}

The disks $D_1, D_2, D_3$ shown in the middle of Figure \ref{family3torosancho3} are properly embedded in $E(K')$. Note that $\bdry D_i = K'$, and $D_i$ intersect $L$ transversely in three points. Let $H_i'$ be the region in $E(K')$ between $D_i$ and $D_{i+1}$, $i\, mod\, 3$. Note that each $H_i'$ determines a simple tangle $(H_i',J',\tau)$, where $J'$ is a preferred longitude of $K'$, and $\tau$ consists of a product arc and a rational tangle $R(m)$, $R(l)$ or $R(-l)$. Consider the double branched covering of $S^3$ along $L$, which is $S^3$ since $L$ is the trivial knot. The lifting of $K'$ under the covering is a knot $K=K(l,m)$ as the linking number between $K'$ and $L$ is an odd number. Observe the following properties: 
\begin{itemize}
    \item The lifting $T_i$ of $D_i$ is a genus one Seifert surface for $K$.
    \item  The lifting $H_i$ of  $H_i'$ is a simple handlebody pair, by Lemma \ref{branchedcover}.
\end{itemize}
    
Let $E_{i,+}, E_{i,+}'$, $E_{i,-}$, $E_{i,-}'$ be disks in $H'_i$ as in Remark \ref{discosenH}: 
   
     The liftings $C_{i,\pm}$ and $C_{i,\pm}'$ of $E_{i,\pm}$ and $E_{i,\pm}'$ are annuli. 
        The core of $C_{i,\pm}'$ is a power curve for $H_i$. 
        Note that both $C_{i,+}$ and $C_{i+1,-}'$, lie in $T_i$, and their cores intersect in one point.  Therefore, the Seifert surfaces $T_i$ determine a decomposition of $E(K)$ into simple pairs that satisfy the hypothesis of Lemma \ref{regionsimple3}.

\begin{theorem}\label{3toricase1}
Let $K$ be a knot of the family of knots $K(l,m)$. Then:  
\begin{enumerate}
     
    \item $E(K)$ has exactly $3$ disjoint, non-isotopic genus one Seifert surfaces.
    \item The Seifert surfaces determine a decomposition of $E(K)$ into simple pairs.
     \item $K$ is a hyperbolic knot.
     \item The circular width of $E(K)$ is $(3)$.
\end{enumerate}
\end{theorem}

\begin{proof}

By the construction of the family $K(l,m)$, the Seifert surfaces $T_1, T_2$ and $T_3$ are genus one and non-isotopic. Two consecutive Seifert tori bound a region $H_i$ which is the lifting of a simple tangle, thus is a simple pair by Lemma \ref{branchedcover}. Moreover, for two adjacent pairs the corresponding power curves intersect in one point. Therefore the conditions of Theorem \ref{theo:MaximalKakimizuSimpleCircularwidth} are satisfied and the items $(1)$-$(3)$ are true.

Claim: Item $(3)$ of Lemma \ref{regionsimple3} is satisfied.

Note that the exterior of $K$, $E(K)$ double branch covers the exterior of $K'$, $E(K')$. There is an involution $\tau$ in $E(K)$, whose rotation axis does not intersect $\bdry E(K)$, such that $E(K)/\tau \cong E(K')$. Perform $0$-Dehn surgery on $K$. This is obtained by filling $\bdry E(K)$ with a solid torus $V$, such that a meridian of $V$ is glued to a preferred longitude of $K$. The involution $\tau$ on $E(K)$ can be extended to an involution $\tau'$ in $E(K)\cup V$, in such a way that the core curve of $V$ remains fixed by $\tau'$. The quotient of this involution corresponds to performing $0$ -Dehn surgery on $K'$. The result is $S^1\times S^2$, with an extra circle in the branch set, as in the right hand side of Figure \ref{family3torosancho3}. We want to apply Lemma \ref{regionsimple3}, with $T_1$ as $T_1$ and $T_4$. Let $H_i$ the simple region bounded by $T_i$ and $T_{i+1}$, for $i=1,2,3$.

The surfaces $T_i$ fill up to surfaces $\hat T_i$ in $E(K)\cup V$. Consider $H_2[J]$ which is just the region between $\hat T_2$ and $\hat T_3$. Note that $H_2[J]$ double branch covers $H_2'[J']$, and this is the tangle in $S^2\times I$ shown in the right side of Figure \ref{family3torosancho3}. 
It is not difficult to check that its double branched cover, i.e., $H_2[J]$, is the torus link $T(2,2m)$, which is also a 2-bridge link. Let $\gamma_1$ and $\gamma_3$ be the power curves of $H_1$ and $H_3$ lying in $T_2$ and $T_3$, respectively. We will show that $\gamma_1$ and $\gamma_3$ are meridian curves of $H_2[J]$. Dehn filling $H_2[J]$ along $\gamma_1$ and $\gamma_3$, correspond to fill up with rational tangles the region $H_2´[J']$. The corresponding filling corresponds to fill up with a tangle such that the essential disk of the tangle is attached along $\bdry E_{i,+}$ and $\bdry E_{i,-}$. The result is shown in Figure \ref{fillingtrivial}. Note that the resulting link is the trivial knot. This implies that the filling of $H_2[J]$ is $S^3$, and then $\gamma_1$ and $\gamma_3$ are meridians, as desired.

Therefore the circular width for this family of knots is $(3)$.
\end{proof}

\subsection{A family of knots with three disjoint non-isotopic Seifert surfaces and circular width (3,3)}

\begin{figure}[h]
\includegraphics[scale=0.8]{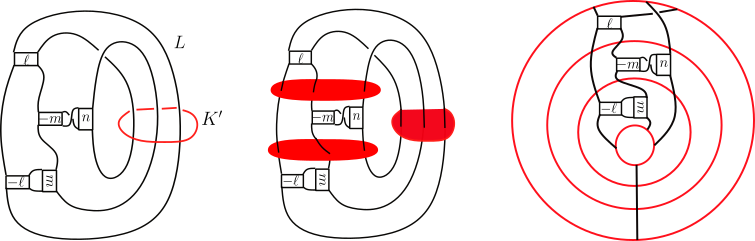}
\caption{A family of knots with three non-isotopic Seifert surfaces and circular width (3,3) 
}
    \label{fig:family3-33}
\end{figure}

Let $K'\cup L$ be the link shown on the left hand side of Figure \ref{fig:family3-33},  where the boxes represent a sequence of  $l$, $m$ or $n$ crossings. Both $K'$ and $L$ are trivial knots.   

\begin{figure}
    \centering
    \includegraphics[width=0.5\linewidth]{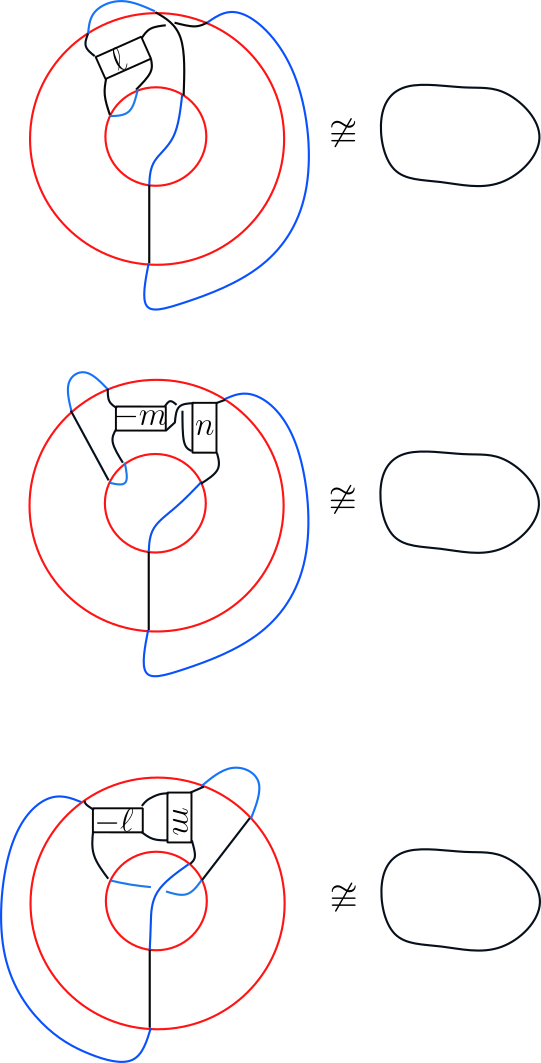}
    \caption{A non-trivial knot}
    \label{non-trivialknots}
\end{figure}

The disks $D_1, D_2, D_3$ 
shown in the middle of Figure \ref{fig:family3-33} are properly embedded in $E(K')$. Note that $\bdry D_i = K'$ and $D_i$ intersect $L$ transversely in three points. Consider the double branched covering of $S^3$ along $L$, which is $S^3$ for $L$ is the trivial knot. 

The lifting of $K'$ under the covering is a knot $K=K(l,m, n)$ as the linking number between $K'$ and $L$ is an odd number. As in the previous example, we have the following properties: 

\begin{itemize}
    \item The lifting $T_i$ of $D_i$ is a genus one Seifert surface for $K$.
    \item  Let $H_i'$ be the region bounded by $D_i$ and $D_{i+1}$. The lifting $H_i$ of  $H_i'$ in $E(K')$ is a genus two handlebody.
    \item Let $E_i, E_i'$ be disks in $H'_i$ such that: 
    \begin{itemize}
        \item $E_i'\subset D_i$ and $E_i' \cap L=$ two points.
        \item $E_i$ is property embedded in $H'_i$ and $E_i \cap L=$ two points.
        \item $\bdry E'_i =\bdry E_i$ 
      
    \end{itemize}
Let $E_{i,+}, E_{i,+}'$, $E_{i,-}$, $E_{i,-}'$ be disks in $H'_i$ as in Remark \ref{discosenH}.
The liftings $C_{i,\pm}$ and $C_{i,\pm}'$ of $E_{i,\pm}$ and $E_{i,\pm}'$ are annuli. 
The core of $C_{i,\pm}'$ is a power curve for $H_i$.  
Note that both $C_{i,+}$ and $C_{i+1,-}'$, lie in $T_i$, and their cores intersect in one point.  Therefore, the Seifert surfaces $T_i$ determine a decomposition of $E(K)$ into simple pairs that satisfy the hypothesis of Lemma \ref{regionsimple3}.

\end{itemize}

\begin{theorem}\label{3toricase2}
Let $K$ be a knot of  the family of knots $K(l,m,n)$. Then:  
\begin{enumerate}

    \item $E(K)$ has exactly $3$ disjoint, non-isotopic genus one Seifert surfaces.
    \item The Seifert surfaces determine a decomposition of $E(K)$ into simple pairs.
     \item $K$ is a hyperbolic knot.
     \item The circular width of $E(K)$ is $(3,3)$.

\end{enumerate}
\end{theorem}
\begin{proof}
By the construction of the family $K(l,m,n)$, the Seifert surfaces $T_1, T_2$ and $T_3$ are genus one and non-isotopic. The proof of items (1), (2) and (3) is identical to the corresponding proof in Theorem \ref{3toricase1}.

Claim: Item $(4)$ of Lemma \ref{regionsimple3} is not satisfied.

We proceed as in the proof of Theorem \ref{3toricase1}. The exterior of $K$, $E(K)$ double branch covers the exterior of $K'$, $E(K')$. There is an involution $\tau$ in $E(K)$ such that $E(K)/\tau \cong E(K')$. Perform $0$-Dehn surgery on $K$. This is obtained by filling $\bdry E(K)$ with a solid torus $V$, such that a meridian of $V$ is glued to a preferred longitude of $K$. The involution $\tau$ on $E(K)$ is extended to an involution $\tau'$ in $E(K)\cup V$, so that that the core curve of $V$ remains fixed by $\tau'$. The quotient of this involution corresponds to performing $0$ -Dehn surgery on $K'$. The result is $S^1\times S^2$, with an extra circle in the branch set, as in the right hand side of Figure \ref{family3torosancho3}. We want to apply Lemma \ref{regionsimple3}, with $T_3$ as $T_1$ and $T_4$. Let $H_i$ the simple region bounded by $T_i$ and $T_{i+1}$, for $i=1,2,3$.

The surfaces $T_i$ fill up to surfaces $\hat T_i$ in $E(K)\cup V$. Consider $H_2[J]$ which is just the region between $\hat T_2$ and $\hat T_3$. Let $\gamma_1$ and $\gamma_3$ be the power curves of $H_1$ and $H_3$ lying in $T_2$ and $T_3$, respectively. We will show that $\gamma_1$ and $\gamma_3$ cannot be meridian curves of $H_2[J]$. Dehn filling $H_2[J]$ along $\gamma_1$ and $\gamma_3$, correspond to fill up with rational tangles the region $H_2'[J']$. The corresponding filling corresponds to filling up with a tangle such that the essential disk of the tangle is attached along $\bdry E_{i,+}$ and $\bdry E_{i,-}$. The result is shown in the top of Figure \ref{non-trivialknots}. Note that the resulting knot is non-trivial. This implies that the filling of $H_3[J]$ is not  $\Sph^3$, and then $\gamma_1$ and $\gamma_3$ are not meridians.
Now repeat the construction, but now taking $T_1$ or $T_2$ as $T_1$. In the middle and bottom parts of Figure \ref{non-trivialknots}, we see that the knot obtained is not trivial. So for any choice of middle region, the knot obtained is non-trivial, and the condition is not satisfied.

Since the simple pairs are as in Lemma \ref{regionsimple}, then the circular width is $(3,3)$.

\end{proof}

\subsection{A family of knots with five disjoint non-isotopic Seifert surfaces of genus one}\label{5toros}

\begin{figure}
    \centering\def\svgwidth{1\linewidth}
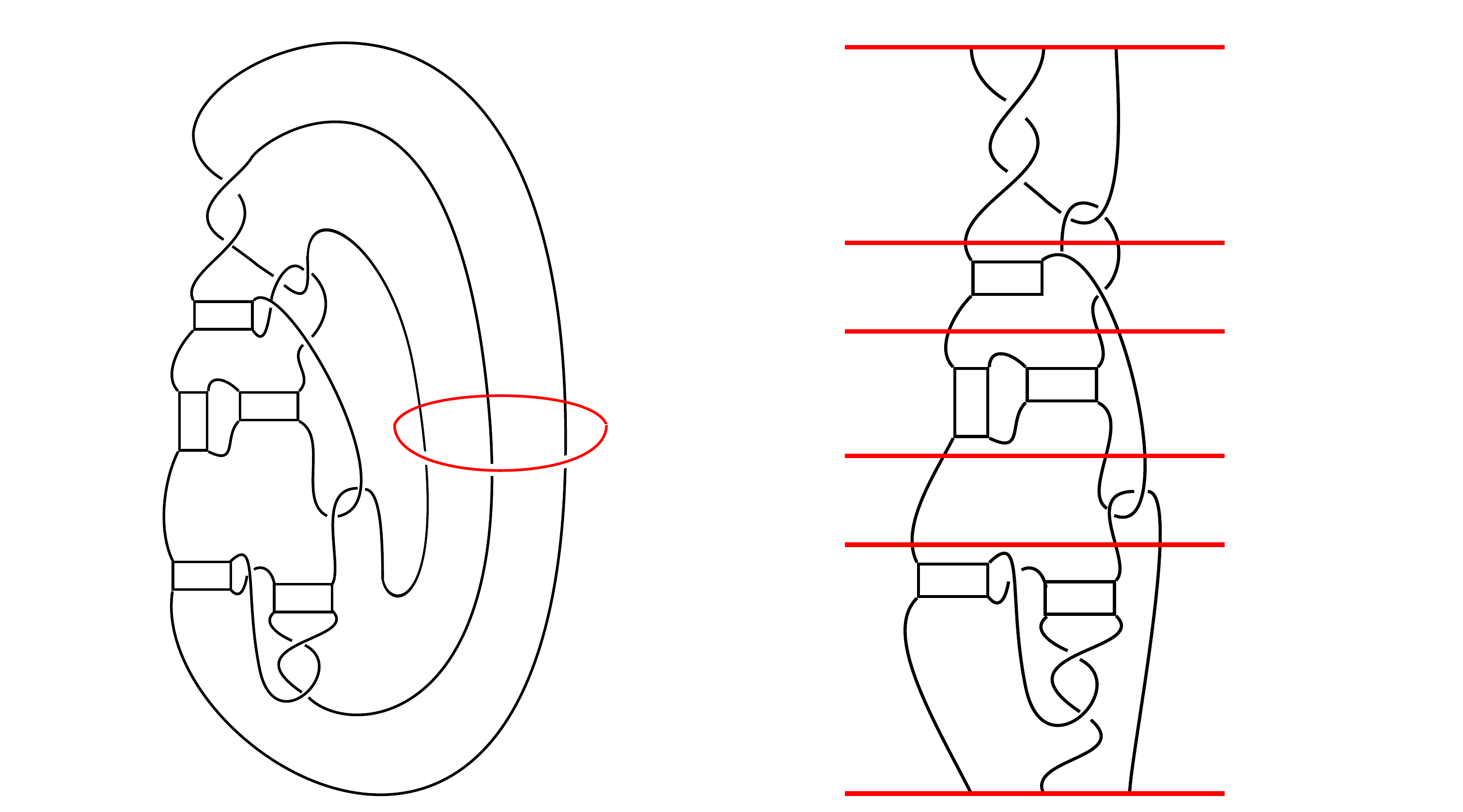
    \caption{A family of knots with five non-isotopic Seifert surfaces}
    \label{family5toros}
\end{figure}

Let $K'\cup L$  be the link shown on the left hand side of Figure \ref{family5toros} where the boxes represent a sequence of  $l$, $m$ or $n$ horizontal crossings. Let us assume that $\vert l\vert, \vert m\vert, \vert n \vert \geq 2$. Note that both $K'$ and $L$ are trivial knots. This link is quite special and not easy to find; it is somehow related to the tangles that determines the family of knots with half integral toroidal Dehn surgery \cite{Mario}.

There are disks $D_1, D_2, D_3, D_4, D_5$, as shown on the right hand side of Figure \ref{family5toros}, which are properly embedded in $E(K')$. Note that $\bdry D_i = K'$ and $D_i$ intersect $L$ transversely in three points. These disks divide $E(K')$ into five regions $H_1'$, $H_2'$, $H_3'$, $H_4'$ and $H_5'$. Each of these regions determines a simple tangle $(H_i',J',\tau)$, which is non-trivial, as follows from the right hand side of Figure \ref{family5toros}. Consider the double branched covering of $S^3$ along $L$, which is $S^3$ since $L$ is the trivial knot. The knot $K=K(l,m,n)$ is the lifting of $K'$ under such covering. The lifting $T_i$ of $D_i$ is a genus one Seifert surface for $K$. The lifting of $H_i'$ is a simple handlebody pair $(H_i,J)$.

\begin{theorem}
Let $K$ be a knot of  the family of knots $K(l,m,n)$. Then:
\begin{enumerate}
    \item  $E(K)$ has exactly $5$ disjoint, non-isotopic genus one Seifert surfaces.
    \item The Seifert surfaces determine a decomposition of $E(K)$ into simple pairs.
    \item $K$ is a hyperbolic knot.
\item The circular width of $E(K)$ satisfies $ cw(E(K)) = (3,3)$.
\end{enumerate}
\end{theorem}
\begin{proof}
By the construction of the family $K(l,m,n)$, the Seifert surfaces $T_1, T_2,...,T_5$ are genus one and non-isotopic. Two consecutive Seifert tori bound a region $H_i$ which is the lifting of a simple tangle, thus is a simple pair by Lemma \ref{branchedcover}. Moreover, for two adjacent pairs the corresponding power curves intersect in one point. Therefore the conditions of Theorem \ref{theo:MaximalKakimizuSimpleCircularwidth} are satisfied and items $(1)$-$(3)$ follow.

We take the decomposition given by the surfaces $T_1$ and $T_4$. This decompose $E(K)$ in two pieces, $M_1=H_1\cup H_2 \cup H_3$ and $M_2=H_4\cup H_5$. By Lemma \ref{regionsimple3} $cw(M_2)=(3)$. Note the middle part of $H_2$, is the double brached cover of a tangle $(H_2',J',\tau)$ given by a rational tangle $R(2)$. By the proof of Theorem \ref{3toricase1} and by Lemma \ref{regionsimple3}, it follows that $cw(M_1)=(3)$. This implies that $cw(E(K))=(3,3)$.
\end{proof}

There is another family of knots $K(l,m,p)$ also coming from the tangles defined in \cite{Mario}, which also have five disjoint, non-isotopic genus one Seifert surfaces.
Presumably, those families of knots coincide with the family of knots constructed by Valdez-Sánchez in \cite{LVgenusonehypknotatmostfive}. We conjecture that those families consist of all hyperbolic knots that have five disjoint, non-isotopic genus one Seifert surfaces. If this conjecture were true, it would imply that the circular width of a hyperbolic genus one knot admitting five genus one Seifert surfaces has to be $(3,3)$. 

\section*{Declaration of IA use}
The authors declare that IA was not used in the generation of research, preparation of statements, proofs and drawings of the present article.  The authors take full responsibility of the mathematical content and proofs in this paper.

\bibliographystyle{alpha}
\bibliography{biblio}

\end{document}